\documentclass[a4paper, 10pt]{amsart}
\usepackage[left=2.5cm, right=2.5cm, 
top=2.5cm, bottom=2.5cm]{geometry}
\usepackage[colorlinks]{hyperref}
\usepackage[english]{babel}

\usepackage{array}
\usepackage{amsmath}
\usepackage{amssymb}
\usepackage{amsthm}
\usepackage{arydshln}
\usepackage{caption}
\usepackage{comment}
\usepackage{diagbox}
\usepackage{enumerate}
\usepackage{float}
\usepackage{graphicx}
\usepackage{mathrsfs}
\usepackage{multirow}
\usepackage{setspace}
\usepackage{stmaryrd}
\usepackage{subfig}
\usepackage{subcaption}
\usepackage{wrapfig}
\usepackage{ulem}
\usepackage{xcolor}
\allowdisplaybreaks[4] 

\numberwithin{equation}{section}
\newcommand{\lj}{[\hspace{-2pt}[}
\newcommand{\rj}{]\hspace{-2pt}]}
\newcommand{\jump}[1]{\lj #1 \rj}
\newcommand{\aver}[1]{ \{#1\}  }
\newcommand{\ds}{\displaystyle}
\newcommand{\ep}{\epsilon}

\newtheorem{lem}{Lemma}
\newtheorem{rmk}{Remark}
\newtheorem{exe}{Example}
\newtheorem{thm}{Theorem}

\newcommand{\bn}{\mathbf{n}}
\newcommand{\pn}{\partial_\mathbf{n}}

\usepackage{etoolbox}
\BeforeBeginEnvironment{comment}{\vspace{-5pt}}
\BeforeBeginEnvironment{figure}{\vspace{-10pt}}
\AfterEndEnvironment{figure}{\vspace{-5pt}}
\BeforeBeginEnvironment{table}{\vspace{-5pt}}
\AfterEndEnvironment{table}{\vspace{-5pt}}
\BeforeBeginEnvironment{center}{\vspace{-5pt}}
\allowdisplaybreaks[4]
\title[Biharmonic Dirichlet Control]
{A Discontinuous Galerkin Method for Biharmonic\\Dirichlet 
Control Problems With Patch Reconstruction}

\author[W. Gong]{Wei Gong}\address{SKLMS \& NCMIS, Institute of 
	Computational Mathematics, Academy of Mathematics and Systems Science, 
	Chinese Academy of Sciences, Beijing 100190, P.R. China}
\email{wgong@lsec.cc.ac.cn}

\author[R. Li]{Ruo Li} \address{CAPT, LMAM and School of Mathematical
	Sciences, Peking University, Beijing 100871, P.R. China;
	Chongqing Research Institute of Big Data, Peking University, Chongqing
	401121, P.R. China}
\email{rli@math.pku.edu.cn}

\author[J. Yin]{Jun Yin} \address{School of Mathematical
	Sciences, Peking University, Beijing 100871, P.R. China}
\email{2401110056@stu.pku.edu.cn}
\date{\today}

\begin{document}
\setlength{\abovedisplayskip}{4pt}
\setlength{\belowdisplayskip}{4pt}
\maketitle

\begin{abstract}
	In this paper, we present and analyze an interior penalty discontinuous Galerkin method for 
	Dirichlet boundary control problems governed by the 
	biharmonic equation, with the control imposed through the normal derivative. The method is based on a reconstructed discontinuous	approximation that admits an arbitrarily high-order approximation space with only one unknown per element.
	We construct a discretization that approximates 
	the state and adjoint state in this space and 
	derive a priori 
	error estimates in the energy and $L^2$ norms.
	The analysis is carried out under the limited regularity available on polygonal domains by means of a medius argument.
	Numerical experiments are provided to support the predicted 	convergence rates.
	
	\noindent \textbf{keywords}: Dirichlet boundary control; 
	biharmonic equation; reconstructed discontinuous 
	approximation
	
\end{abstract}

\section{Introduction}\label{sec1}
The numerical approximation of optimal control problems 
constrained by partial differential equations (PDEs) represents a cornerstone of modern computational science, driven by diverse applications ranging from flow control and shape optimization to aerospace engineering. These problems seek to determine a control variable that minimizes a given  cost functional while satisfying a PDE state equation. Depending on the nature of the actuation, they are broadly classified into distributed control problems, where the control acts as an external force within the domain, and boundary control problems, where the control enters through Dirichlet, Neumann, or Robin boundary conditions. The mathematical theory of optimal control problems is presented in \cite{lions1971ocpde}. The numerical approximation of this class of problems varies depending 
on the type of control; for a comprehensive overview 
of numerical methods, we refer to 
\cite{hinze2008ocpde,manzoni2021ocpde}. Among the 
various numerical strategies, finite element methods (FEMs) 
have emerged as the most compelling approach due to their 
efficiency, flexibility, and natural adaptivity for complex 
geometries. 
The foundational work on finite element 
approximation and rigorous error analysis for this class of 
problems began with the work of Falk \cite{falk1973ocp}
and Geveci \cite{geveci1979ocp} in the 1970s, 
a paradigm that has since inspired extensive research and 
refinement, as documented in the extensive literature cited 
therein. 
More specifically, for distributed controls, we refer to \cite{aradacasastroltzsch,  
	frei2013bidistributed,hinze2005varydiscre, li2007posteriori}.
For the Neumann boundary control, we refer to \cite{brenner2019neumann,  
casas2012neumann,casas2008neumann, geveci1979ocp}. 
For Dirichlet boundary control, we refer to \cite{casas2006dirichlet,chowdhury2017dirichlet,gong2022dirichlet,may2013dirichlet,of2015dirichlet} and the references therein.
	
Although optimal control 
problems with second-order state equation 
have been widely investigated 
\cite{brenner2019neumann,casas2012neumann,casas2008neumann,chowdhury2017dirichlet,falk1973ocp, geveci1979ocp, li2007posteriori}, 
	there is relatively less literature that concentrates on the analysis of those whose governing equation is of higher order. The first result on optimal control of biharmonic equations appeared in \cite{caoyangbiharmonic}, where the authors used the Ciarlet-Raviart mixed scheme and derived a priori error estimates. 
In \cite{chowdhury2014bidistributed},
a $C^0$-interior penalty method is analyzed for 
distributed control 
governed by a biharmonic equation with box-constrained control 
on general polygonal domains. 
An optimal order a priori error estimate has been derived 
with a superconvergence estimate. Furthermore, a
residual based a posteriori error estimator is derived therein, 
which is
proved to be reliable and efficient. The 
authors in \cite{frei2013bidistributed} have studied the same  
problem with a mixed formulation in convex polygonal domains.
They have proposed a mixed finite element method and derived 
the optimal order a priori
error estimate through a superconvergence result 
of the optimal control as well. 
In \cite{chowdhury2015priori}, the biharmonic
problem with box-constrained control is investigated.
A general framework for 
a priori and a posteriori error estimates of discontinuous 
Galerkin (DG) methods is derived for 
both distributed and Neumann control on general polygonal domains. 

The analysis of Dirichlet boundary control is 
comparatively more 
subtle than that of distributed and Neumann boundary control,
because the control does not enter the variational formulation in a natural way. 
Several approaches have been proposed to overcome 
this difficulty. 
One possibility is to seek the control in $L^2(
\Gamma)$ as in \cite{casas2006dirichlet, 
may2013dirichlet}, which leads to a
very weak (transposition) solution. It is then difficult to derive optimal-order error estimates, especially in the $L^2$ norm. 
Alternatively, \cite{of2015dirichlet,gong2022dirichlet} penalize the control in 
$H^{\frac12}(\Gamma)$, where the Steklov--Poincar\'e operator 
is used to realize the Dirichlet-to-Neumann map via a 
harmonic extension. Error estimates are then derived in 
the energy and $L^2$ norms. 
The authors of 
\cite{chowdhury2017dirichlet, gong2019adaptive} penalize the control in $H^1(\Omega)$, which is equivalent to $H^{\frac12}(\Gamma)$, 
and derive optimal-order a priori error estimates 
in the energy and $L^2$ norms, 
as well as an a posteriori error analysis. 
For problems governed by biharmonic equations, 
\cite{chowdhury2017bipriori, chowdhury2024biprioripos} 
seek the control in $H^2(\Omega)$ and develop a $C^0$ interior 
penalty method, together with a priori and a posteriori 
error estimates and 
reliable and efficient error 
estimators.
Recently, in \cite{garg2025biprioripos}, 
the authors proposed a discontinuous Galerkin 
method for the 
same problem with a corresponding error analysis. 

Notice that the Dirichlet boundary controls in 
\cite{chowdhury2017bipriori, chowdhury2024biprioripos,
	garg2025biprioripos} 
are all imposed on the Dirichlet trace.
In this article, we consider a related problem with
the same governing equation, 
but with the control acting on the  
normal derivative. Specifically, we consider
\begin{equation}\label{eq2}
	\begin{aligned}
		&\min_{u\in U_{ad}}\frac{1}{2}\|\psi-\psi_d\|_{
			L^2(\Omega)}^2+\frac\alpha2
		\|u\|_U^2,&\\
		\mathrm{s.t.}\ 
		&\Delta^2 \psi= f\quad\mathrm{in}\ 
		\Omega,\\
		&\psi=0,\ \pn\psi 
		=u\quad \mathrm{on}\ \Gamma:=\partial\Omega,
	\end{aligned}
\end{equation}
where $u$ is the control, $\psi$ is the state,  $U$ is a control space and $U_{ad}\subseteq U$ is a subspace that will be specified later. Further details on \eqref{eq2}
are given in Section \ref{sec2}. 
We develop an a priori error analysis for this problem, deriving energy-norm error estimates
on any convex polygonal domain and 
$L^2$ norm estimates when
the interior angles are less than 126.283 degrees. 

Although Dirichlet boundary control of the  biharmonic equation \eqref{eq2} is of independent interest, another important motivation comes from the study of Dirichlet boundary control for Stokes flows. 
Consider the following optimal control problem:
\begin{equation}\label{eq1}
\begin{aligned}
&\min_{u\in U_{ad}}\frac{1}{2}\|{\bf y}-{\bf y}_d\|_{L^2(\Omega)}^2+\frac\alpha2
\|u\|_U^2,&\\
\mathrm{s.t.}\ &-\nu\Delta {\bf y}+\nabla p={\bf f},\ \nabla\cdot {\bf y}=0\quad\mathrm{in}\ 
\Omega,\\
&{\bf y}=u\boldsymbol{\tau}\quad\mathrm{on}\ \Gamma,
\end{aligned}
\end{equation}
where $\nu$ is a positive constant, ${\bf y}$ is the velocity and $p$ is the pressure, $u$ is the scalar boundary control, ${\bf f}$ and ${\bf y}_d$ are given functions. Here we consider only tangential boundary controls where the normal velocity vanishes on the boundary (cf. \cite{chen2023dirichlet,hu2020dirichlet}). Let $\psi$ be a stream function such that ${\bf y}=(\frac{\partial\psi}
{\partial x_2},-\frac
{\partial\psi}{\partial x_1})^T.$ Then the control problem 
becomes
\begin{equation}\label{eq22}
\begin{aligned}
&\min_{u\in U_{ad}}\frac{1}{2}\|\nabla\psi-{\bf z}_d\|_{L^2(\Omega)}^2+\frac\alpha2
\|u\|_U^2,&\\
\mathrm{s.t.}\ 
&\nu\Delta^2 \psi=\nabla\times {\bf f}\quad\mathrm{in}\ 
\Omega,\\
&\psi=0,\ \pn\psi 
=-u\quad \mathrm{on}\ \Gamma,
\end{aligned}
\end{equation}
where ${\bf z}_d:=(-{\bf y}_d^2,{\bf y}_d^1)^T$. 
With this stream-function approach, the divergence-free constraint is automatically satisfied at the discrete level, in contrast to \cite{chen2023dirichlet}, where a special scheme was designed for that purpose. 


To solve the problem \eqref{eq2} numerically, we adopt an efficient DG 
method called reconstructed discontinuous approximation (RDA) 
proposed in \cite{li2024RDA,li2024biharm}. On the one hand, 
conforming finite elements for fourth order problems 
force the finite element space to be contained in 
$H^2(\Omega)$, whose implementation is rather inefficient. 
The DG methods are inherently well-suited for high order 
partial differential equations, owing to their intrinsic 
properties of local conservation and excellent parallel 
scalability. On the other hand, 
to achieve high order accuracy,
traditional DG methods are much more computationally 
expensive
compared to continuous elements because of the 
significant increase in the number of degrees
of freedom. The RDA method in \cite{li2024RDA} 
allows us to construct a high order
approximation space requiring only one unknown per element. 
The resulting space forms a small subspace
of the standard discontinuous piecewise polynomial space, 
which is linearly isomorphic to the piecewise 
constant space. Therefore, the number of degrees of 
freedom is reduced to the number of elements in the mesh. 
Nevertheless, the approximation property of this subspace 
has been proved to be optimal under certain conditions.

The rest of this article is organized as follows. In Section \ref{sec2}, we present the model problem,  investigate the well-posedness of the optimization problem and derive the regularity of solutions. In Section \ref{sec_numer} we describe
the procedure and features of the RDA method, and introduce the discontinuous Galerkin approximation to biharmonic equations. We also put 
forward the discrete first order optimality system to be solved numerically. In Section 
\ref{sec3}, we prove the a priori error
estimates under the energy and $L^2$ norms.
In Section \ref{sec4}, we present some numerical 
examples to support the theoretical results.
Finally, we summarize the article in Section \ref{sec5}.

\section{Optimal Control Problems}\label{sec2}

Let $\Omega\subseteq
\mathbb{R}^2$ be a convex polygonal domain with boundary $\Gamma=
\partial\Omega$ consisting of $l$ segments $\Gamma_1,\Gamma_2,
\cdots,\Gamma_l.$ $U$ is a properly chosen control space. 
$\alpha>0$ is the regularization parameter. 
$f,\psi_d\in L^2(\Omega)$ are given functions. 
The state variable 
$\psi = \psi(u,f)$ depends on $u$ and $f$ uniquely, so the cost 
functional is essentially a quadratic functional with respect to $u$ since the governing equation is linear.
%

The choice of the control space determines the formulation of
the state equation and the regularity of the optimal solution.
We consider $U=L^2(\Gamma)$ and $U=H^1(\Gamma)$ in a unified
framework. The energy-space choice $U=H^{\frac12}(\Gamma)$ will be
treated elsewhere. In the rest of this article,
we use $(\,\cdot\,,\,\cdot\,)$ to 
denote the $L^2$ inner 
product in $\Omega$ and $(\,\cdot\,,\,\cdot\,)_\Gamma$ to 
denote the $L^2$ inner product on $\Gamma$.
Let $\mathcal V_\Gamma$ denote the set of vertices of the polygon.
The admissible space is a closed linear subspace
\begin{equation}\label{eq:Uad}
U_{ad}:=
\begin{cases}
L^2(\Gamma),&U=L^2(\Gamma),\\
\{v\in H^1(\Gamma)\mid v(P)=0,\forall P\in\mathcal V_\Gamma\},
&U=H^1(\Gamma).
\end{cases}
\end{equation}
Here $H^1(\Gamma)$ is the intrinsic Sobolev space on the closed
polygonal curve. Its edgewise $H^1$ functions have matching
values at each corner, and $H^1(\Gamma)\subseteq C(\Gamma)$.
We also write
\[
H_{\rm pw}^s(\Gamma):=\prod_{i=1}^lH^s(\Gamma_i).
\]
Within the class $u\in H^1(\Gamma)$, the vertex conditions in
\eqref{eq:Uad} imply the polygonal compatibility conditions for the
pair $(\psi,\pn\psi)=(0,u)$ to be the trace of a function in
$H^2(\Omega)$.
Define a bilinear form 
\[a(\psi,\phi):=(\Delta\psi,\Delta\phi),
\quad\forall\psi,\phi\in H^2(\Omega).\]
We will detail these two cases respectively.

\begin{itemize}
    \item $U:=L^2(\Gamma)$. Since we are concerned with 
    biharmonic equations, the standard energy space $\psi\in 
    H^2(\Omega)$ requires a Dirichlet boundary condition 
    $u\in H^{1\over 2}(\Gamma)$, so that the governing state 
    equation has to be understood in the very weak sense. 
    That is, we are looking for $\psi\in H^1_0(\Omega)$ such that
\begin{equation}\label{veryweaksolution}
      (\nabla\psi,\nabla\Delta\phi)=-(f,\phi)+( u,\Delta \phi)_\Gamma, 
      \quad\forall \phi\in H^3(\Omega)\cap H_0^2(\Omega).
    \end{equation}
    This is the transposition solution for rough boundary data.
    Its well-posedness follows from the  solution operator of biharmonic equations with clamped boundary conditions.
    For $G\in H^{-1}(\Omega)$, let $\phi_G\in
    H^3(\Omega)\cap H_0^2(\Omega)$ solve the following problem
    \[
    \Delta^2\phi_G=G,\quad
    \|\phi_G\|_{H^3(\Omega)}\leq C\|G\|_{H^{-1}(\Omega)}.
    \]
    Define $\psi\in H_0^1(\Omega)$ by
    \[
    \langle G,\psi\rangle_{H^{-1},H_0^1}
    :=(f,\phi_G)-(u,\Delta\phi_G)_\Gamma .
    \]
    The trace theorem and the above regularity estimate show that
    the right-hand side is a bounded linear functional on $H^{-1}(\Omega)$.
    Since $(H^{-1}(\Omega))'=H_0^1(\Omega)$, this defines a unique
    $\psi$ and gives
    \[
    \|\psi\|_{H^1(\Omega)}
    \leq C\bigl(\|f\|_{L^2(\Omega)}
    +\|u\|_{L^2(\Gamma)}\bigr).
    \]
    Taking $G=\Delta^2\phi$ recovers \eqref{veryweaksolution}.
    Now we can introduce the affine control to state mapping $S:L^2(\Gamma)\rightarrow H^1_0(\Omega)$ by 
setting $\psi:=Su$, and we obtain the reduced optimization problem
\begin{equation}
    \min_{u\in U_{ad}}\ \hat{J}(u):=\frac{1}{2}\|Su-\psi_d
    \|_{L^2(\Omega)}^2+\frac\alpha2
\|u\|_{L^2(\Gamma)}^2.
\end{equation}
Then, using standard arguments, we can prove the existence of  a unique solution to the optimal control problem.

\item $U:=H^{1}(\Gamma)$. For every $u\in U_{ad}$, the corner
compatibility in \eqref{eq:Uad} provides an $H^2$ lifting of the boundary data 
$(0,u)$, and hence the state equation has a standard variational
solution $\psi\in H^2(\Omega)$.
 That is, find $\psi\in H^2(\Omega)$ such that $\psi=0$ and $\pn\psi=u$ on $\Gamma$ and 
    \begin{equation}
        a(\psi,\phi)=(f,\phi),\quad\forall \phi\in H_0^2(\Omega).
    \end{equation}
    We can therefore introduce the affine control to state mapping $S:U_{ad}
    \rightarrow H^2(\Omega)$ by setting $\psi:=Su$, and we obtain the 
    reduced optimization problem
\begin{equation}
    \min_{u\in U_{ad}}\ \hat{J}(u):=\frac{1}{2}\|Su-\psi_d\|_{L^2
    	(\Omega)}^2+\frac\alpha2
\|u\|_{H^{1}(\Gamma)}^2.
\end{equation}
Then, by standard arguments, we can also prove the existence of 
a unique solution to the optimal control problem. 
\end{itemize}
 
To derive the optimality condition of problem $\eqref{eq2}$, 
we separate
the state equation into two independent parts as
\begin{flalign*}
\left\{
\begin{aligned}
&\Delta^2\psi_f= f\quad\mathrm{in}\ \Omega,\\
&\psi_f=0,\ \pn\psi_f=0\quad \mathrm{on}\ \Gamma,
\end{aligned}
\right.\quad
\left\{
\begin{aligned}
&\Delta^2\psi_0(u)=0\quad\mathrm{in}\ \Omega,\\
&\psi_0(u)=0,\ 
\pn\psi_0(u)=u\quad \mathrm{on}\ \Gamma.
\end{aligned}
\right.
\end{flalign*}
It is clear that the mapping between $u$ and $\psi_0(u)$ is linear, bounded, and Fr\'echet differentiable, and we denote $\psi_0'(u)v$ its Fr\'echet derivative at $u$ along the direction $v$. Then it is clear that $\psi_0'(u)v=\psi_0(v)$.

Denote the optimal solution to $\eqref{eq2}$ by $(\psi^\ast,
u^\ast)$, then 
\begin{equation}\label{eq4}
\left\{
\begin{aligned}
&\Delta^2\psi^*= f\quad\mathrm{in}\ \Omega,\\
&\psi^*=0,\ 
\pn\psi^*=u^*\quad \mathrm{on}\ \Gamma,
\end{aligned}
\right.
\end{equation}
and $\psi_0^*:=\psi_0(u^*)=\psi^\ast-\psi_f$. 
Define the adjoint variable $\eta^*$ such that 
\begin{equation}\label{eq5}
\left\{
\begin{aligned}
&\Delta^2\eta^*=\psi^*-\psi_d\quad\mathrm{in}\ \Omega,\\
&\eta^*=0,\ 
\pn\eta^*=0\quad \mathrm{on}\ \Gamma.
\end{aligned}
\right.
\end{equation}

Denote the cost functional by $J(\psi_0(u),u).$ Then 
$$J'(\psi_0^*,u^*)(v)=(\psi_0^*+\psi_f-\psi_d,
\psi_0'(u^*)v)+\alpha(u^*,v)_U=0,\quad\forall v\in U_{ad}.$$
Notice that $\psi_0'(u^*)(v)=\psi_0(v)$ for any $v\in U_{ad}$.
Therefore,
\[
(\Delta^2\eta^*,\psi_0(v))+\alpha(u^*,v)_U=0,
\quad\forall v\in U_{ad}.
\]
For $U=L^2(\Gamma)$, $\psi_0(v)$ is only a transposition solution in
$H_0^1(\Omega)$, so the classical twice-integrated Green formula is
not available. Instead, \eqref{veryweaksolution} with $f=0$ and
$\phi=\eta^*$ gives
\[
(\nabla\psi_0(v),\nabla\Delta\eta^*)=(v,\Delta\eta^*)_\Gamma.
\]
Distributional integration by parts, using
$\psi_0(v)\in H_0^1(\Omega)$, gives
\[
(\Delta^2\eta^*,\psi_0(v))
=-(\nabla\Delta\eta^*,\nabla\psi_0(v))
=-(\Delta\eta^*,v)_\Gamma.
\]
For $U=H^1(\Gamma)$ the same identity follows from the standard
Green formula because $v\in U_{ad}$ admits an $H^2$ state lifting.
Consequently,
\[
\alpha(u^*,v)_U-(\Delta\eta^*,v)_\Gamma=0,
\quad\forall v\in U_{ad}.
\]
So we have derived the first order optimality system for problem 
$\eqref{eq2}$ as follows:
\begin{equation}\label{eq6}
	\left\{
	\begin{aligned}
		&\Delta^2\psi^*=f\quad\mathrm{in}\ \Omega,\quad
		\psi^*=0,\ \pn
		\psi^*=u^*\quad \mathrm{on}\ \Gamma,\\
		&\Delta^2\eta^*=\psi^*-\psi_d\quad\mathrm{in}\ \Omega,\quad
		\eta^*=0,\ \pn\eta^*=0\quad \mathrm{on}\ \Gamma,\\
		&\alpha(u^*,v)_U-(\Delta\eta^*,
		v)_\Gamma=0,\quad \forall v\in U_{ad}.
	\end{aligned}
	\right.
\end{equation}

In the following theorem we state the regularity of solutions to the optimal control problem \eqref{eq2}.
\begin{thm}\label{thm:regularity}
Assume that $\Omega\subseteq\mathbb{R}^2$ is a bounded convex polygon and $f,\psi_d\in L^2(\Omega)$. Let $(\psi^*,u^*,\eta^*)$ be the optimal solution to the optimal control problem \eqref{eq2}. Then we have the following regularity:
\begin{eqnarray}
    U=L^2(\Gamma):\quad u^*\in H^{1\over 2}(\Gamma),\quad \psi^*\in H^2(\Omega)\cap H_0^1(\Omega),\quad \eta^*\in H^3(\Omega)\cap H_0^2(\Omega);\\
  U=H^1(\Gamma):\quad   u^*\in H_{\rm pw}^{2}(\Gamma)\cap U_{ad},\quad \psi^*\in H^{3-\delta}(\Omega)\cap H_0^1(\Omega),\quad \eta^*\in H^3(\Omega)\cap H_0^2(\Omega).
    \end{eqnarray}
    Furthermore, if all interior angles of $\Omega$ are
less than 126.283 degrees, we have 
\begin{eqnarray}
    U=L^2(\Gamma):\quad u^*\in \gamma_0(H^2(\Omega)),\quad \psi^*\in H^3(\Omega)\cap H_0^1(\Omega),\quad \eta^*\in H^4(\Omega)\cap H_0^2(\Omega);\\
    U=H^1(\Gamma):\quad u^*\in H_{\rm pw}^{2}(\Gamma)\cap U_{ad},\quad \psi^*\in H^{3-\delta}(\Omega)\cap H_0^1(\Omega),\quad \eta^*\in H^4(\Omega)\cap H_0^2(\Omega),
    \end{eqnarray}
    where $\gamma_0(\cdot)$ is the Dirichlet trace operator and $\delta>0$ is arbitrary.
    \end{thm}
\begin{proof}
We use the following regularity estimates for biharmonic equations with clamped boundary conditions on a convex polygon
\cite{blumrannacher,grisvard2011regupara}: the solution operator maps
$H^{-2}(\Omega)$ to $H_0^2(\Omega)$ and $H^{-1}(\Omega)$ to
$H^3(\Omega)\cap H_0^2(\Omega)$ continuously. Under the stated
interior angle restriction it also maps $L^2(\Omega)$ to
$H^4(\Omega)\cap H_0^2(\Omega)$ continuously.

We first construct a lifting that respects the polygonal corners.
There is a vector field ${\bf b}\in C^\infty(\overline\Omega)^2$ such that
${\bf b}\cdot\bn=1$ on every open edge. Indeed, at a vertex with
incident normals $\bn_1$ and $\bn_2$, we can choose a constant vector ${\bf b}_P$
solving ${\bf b}_P\cdot\bn_1={\bf b}_P\cdot\bn_2=1$; ${\bf b}_P$ exists uniquely because  the normals are linearly
independent. Away from the vertices we can choose constant vectors with
the required normal component on the single adjacent edge.
A smooth partition of unity subordinate to these boundary
neighborhoods and an interior neighborhood gives $\bf b$.

For $q\in H^{k+1}(\Omega)\cap H_0^2(\Omega)$, $k=2,3$, we set 
$G_q={\bf b}\cdot\nabla q$. On each open edge,  the clamped boundary conditions for $q$ imply
$\nabla q=0$, $\partial_{\tau\tau}q=\partial_{\tau n}q=0$,
and hence $D^2q=(\Delta q)\bn\otimes\bn$. Consequently
\begin{equation}\label{eq:clamped-trace}
\begin{gathered}
G_q\in H^k(\Omega)\cap H_0^1(\Omega),\quad
\pn G_q=\gamma_0(\Delta q),\quad 
\|G_q\|_{H^k(\Omega)}\leq C\|q\|_{H^{k+1}(\Omega)}.
\end{gathered}
\end{equation}
Here, the normal-trace identity follows from
$\pn({\bf b}\cdot\nabla q)=(\pn {\bf b})\cdot\nabla q+{\bf b}\cdot D^2q\,\bn$.
Thus, \eqref{eq:clamped-trace} provides the full corner-compatible
lifting.

Let $U=L^2(\Gamma)$. The definition of very weak solutions to the state equation by transposition methods first gives
$\psi^*\in H_0^1(\Omega)$. Therefore, 
$\psi^*-\psi_d\in L^2(\Omega)$ and
$\eta^*\in H^3(\Omega)\cap H_0^2(\Omega)$ using  standard regularity results. The first order optimality condition yields
$u^*=\alpha^{-1}\gamma_0(\Delta\eta^*)$. In particular
$u^*\in H^{\frac12}(\Gamma)$, since $\Delta\eta^*\in H^1(\Omega)$.
By \eqref{eq:clamped-trace}, $G=\alpha^{-1}G_{\eta^*}$ is an
$H^2$ lifting of the boundary data $(0,u^*)$. The variational problem
\[
a(v,\phi)=(f,\phi)-a(G,\phi),\quad\forall\phi\in H_0^2(\Omega),
\]
has a unique solution $v\in H_0^2(\Omega)$. The function $G+v$
satisfies \eqref{veryweaksolution}, as follows by integrating once
against $\Delta\phi$ for $\phi\in H^3(\Omega)\cap H_0^2(\Omega)$.
The uniqueness of the transposition solution gives $\psi^*=G+v\in H^2(\Omega)$.

Under angle restriction, we can expect $\eta^*\in H^4(\Omega)\cap H_0^2(\Omega)$, from which $u^*\in\gamma_0(H^2(\Omega))$ and the same lifting $G$ belongs
to $H^3(\Omega)$. Since
$f-\Delta^2G\in H^{-1}(\Omega)$, the standard regularity for biharmonic equations with clamped boundary conditions applied to
$v=\psi^*-G$ gives $\psi^*\in H^3(\Omega)\cap H_0^1(\Omega)$.
We note that this argument does not infer critical-order trace compatibility
merely from vanishing corner values.

Now, let $U=H^1(\Gamma)$. The definition of $U_{ad}$ provides an
$H^2$ state lifting, so $\psi^*\in H^2(\Omega)\cap H_0^1(\Omega)$ initially.
The same regularity arguments for biharmonic equations with clamped boundary conditions give $\eta^*\in H^3(\Omega)\cap H_0^2(\Omega)$,
and $\eta^*\in H^4(\Omega)\cap H_0^2(\Omega)$ under angle restriction.
Testing the first order optimality condition with functions supported on each macro-edge
and using their zero endpoint values gives
\[
-\alpha\,\partial_s^2u^*+\alpha u^*=\Delta\eta^*
\quad\text{on }\Gamma_i,\quad
u^*=0\quad\text{at the endpoints of }\Gamma_i.
\]
The right-hand side is in $L^2(\Gamma_i)$. The one-dimensional
Dirichlet regularity therefore gives
$u^*\in H_{\rm pw}^2(\Gamma)\cap U_{ad}$.

For completeness, the polygonal trace theorem
\cite[Chap.~1, Sec.~1.5]{grisvard2011regupara} gives, for
$0<\delta<1/2$ and $v\in H_{\rm pw}^2(\Gamma)\cap U_{ad}$,
a lifting $L_\delta v$ satisfying
\begin{equation}\label{eq:fractional-lifting}
\begin{gathered}
L_\delta v\in H^{3-\delta}(\Omega)\cap H_0^1(\Omega),\quad
\pn L_\delta v=v,\quad 
\|L_\delta v\|_{H^{3-\delta}(\Omega)}
\leq C_\delta\|v\|_{H_{\rm pw}^2(\Gamma)}.
\end{gathered}
\end{equation}
Indeed, the prescribed value trace is zero, the gradient trace
$v\bn$ has zero matching vertex values, and the second-derivative
trace has order $1/2-\delta<1/2$, at which no further
vertex matching is required. Interpolation between the clamped 
$H^{-2}(\Omega)\to H^2(\Omega)$ and $H^{-1}(\Omega)\to H^3(\Omega)$ estimates gives the
$H^{-1-\delta}(\Omega)\to H^{3-\delta}(\Omega)$ shift. Applying this result to
$\psi^*-L_\delta u^*$, whose right-hand side is
$f-\Delta^2L_\delta u^*\in H^{-1-\delta}(\Omega)$, we can obtain
the claimed regularity. The conclusions for larger $\delta>0$
follow by the Sobolev embedding.
\end{proof}

\section{Discontinuous Galerkin Approximation}\label{sec_numer}

\subsection{The Reconstruction Operator}\label{sec2.2}
At first, we give the main idea on the construction of the finite element space.
Let $\mathcal{T}_h$ be a family of shape-regular and quasi-uniform triangulations of 
$\Omega\subseteq\mathbb{R}^2$ with mesh size 
$h$. Let $\mathcal{E}_h^I$ and 
$\mathcal{E}_h^B$ be the set of all  one dimensional 
edges in the interior of $\Omega$ in $\mathcal{T}_h$, and edges 
on the boundary $\Gamma$ in $\mathcal{T}_h$ 
respectively. 
The set of all edges is denoted by 
$\mathcal{E}_h:=\mathcal{E}_h^I\cup\mathcal{E}_h^B.$
Similarly, let $\mathcal{N}_h^I$ and 
$\mathcal{N}_h^B$ be the set of nodes 
in the interior of $\Omega$ in $\mathcal{T}_h$, and nodes 
on the boundary $\Gamma$ in $\mathcal{T}_h$ 
respectively. 
The set of all nodes is denoted by 
$\mathcal{N}_h:=\mathcal{N}_h^I\cup\mathcal{N}_h^B.$
Now we will briefly introduce the reconstruction 
operator $\mathcal{R}^{m}$ put forward in \cite{li2024RDA}. The 
whole process includes three steps in total.

{\it Step 1.} Given $K\in\mathcal{T}_h$, we want to construct an 
element patch $S(K)$ consisting of $K$ itself and some of its
neighbors. Generally, $S(K)$ is 
constructed by continuously adding neighboring elements. 
We first let $S_0(K)=\{K\}$, and then define $S_t(K)$ 
in the following recursive way.
$$S_t(K)=\bigcup\limits_{K'\in S_{t-1}(K)}U(K'),\quad 
t=1,2,\cdots,$$ 
where $$U( K) : = \{ K^{\prime }\in \mathcal{T}_h
\mid\exists e\in\mathcal{E}_h,
e\subseteq\partial{K^{\prime }}\cap 
\partial {K}\} .$$
The process repeats until $t$ reaches the 
condition that $\# S_{t}( K) \geq \# S$, where $\#S$
is a given constant independent of $K$. Then
we can define the patch $S( K)$ as $S_{t}( K) .$ 
By numerical observations, the threshold $\# S$ cannot 
be too small. 
Practically, it can be 
chosen as $1.5\times\mathrm{dim}(\mathbb{P}_m(K))$.

{\it Step 2.} For each $K\in\mathcal{T}_h$, we need to solve a local 
least squares problem in $S(K)$. Let
$\boldsymbol{x}_K$ be the barycenter of $K$ and 
$X(K):=\{\boldsymbol{x}_{K'}\in\Omega\mid K'\in S(K)\}$ be 
the set of barycenters of elements in $S(K)$.  
Let $\Phi_h^{0}$ be the following piecewise constant space  on $\mathcal{T}_h$
$$\Phi_h^0:=\{z_h\in L^2(\Omega)\mid z_h|_K\in\mathbb
{P}_0(K),\:\forall K\in\mathcal{T}_h\}.$$

Given $z_h\in \Phi_h^0$, for each 
$K\in\mathcal{T}_h$, we want to find a polynomial of 
degree $m$ by solving a constrained 
least squares problem:
\begin{equation*}
\begin{aligned}
&\underset{p_m\in\mathbb{P}_m(S(K))}{\operatorname*{
\min}}\sum_{\boldsymbol{x}\in X(K)}(p_m(\boldsymbol{x})
-z_h(\boldsymbol{x}))^2,\\
\mathrm{s.t.}\ &p_m(\boldsymbol{x}_K)=z_h(\boldsymbol{x}_K).
\end{aligned}
\end{equation*}
The existence of a solution to this problem is 
induced by the finite dimensionality of $\mathbb{P}_m(S(K)).$
Uniqueness can be proved under certain assumptions 
(cf. \cite[Assumption 1 and Lemma 1]{li2024RDA}).
Denote this solution by $z_{h,S(K)}
\in\mathbb{P}_m(S(K))$. 

{\it Step 3.}  Given $K\in\mathcal{T}_h$, notice that $z_{h,S(K)}$ 
depends linearly  on the function $z_h$, we can then define a 
linear operator $\mathcal{R}_K^m:\Phi_h^0
\longrightarrow\mathbb{P}_m(S(K))$ that maps $z_h$ to $z_{h,S(K)}$.
Finally, we define the global 
reconstruction operator $\mathcal{R}^m$ as
$$\begin{aligned}
&\mathcal{R}^{m}:\Phi_h^0\longrightarrow \Phi_h^m\\
\mathrm{s.t.}\ &(\mathcal{R}^{m}z_h)|_{K}:=(\mathcal{R}_{K}^{m}
z_h)|_{K},\quad\forall K\in\mathcal{T}_h,z_h\in\Phi_h^0.
\end{aligned}$$
Here $\Phi_h^m:=\mathcal{R}^
m\Phi_h^0$ is the image space. $(\mathcal{R}^mz_h)|_K$ is the 
restriction of $\mathcal{R}_K^mz_h$ on $K\in\mathcal{T}_h$. 
Intuitively, $\mathcal{R}^m$ can be regarded as a combination of 
$\mathcal{R}_K^m$, which associates a piecewise constant function 
with a piecewise polynomial of $m$-th order. 

It can be shown 
that $\mathcal{R}^m$ is a linear isomorphism from 
$\Phi_h^0$ to $\Phi_h^m$, i.e., $\mathrm{dim}(\Phi_h^m)=\mathrm{dim}
(\Phi_h^0)$. Moreover, the operator $\mathcal{R}^m$ has an optimal approximation property. For every $s>1$, define $I_h:H^{s}(\Omega)\to \Phi_h^0$ as the 
interpolation such that 
$$I_hz(\boldsymbol{x}_K)=z(\boldsymbol{x}_K),\quad\forall K\in\mathcal{T}_h
,z\in H^{s}(\Omega).$$ 
Throughout the analysis $m$ is fixed and the reconstruction
patches satisfy uniform approximation hypotheses: their diameters
are comparable to $h$, they have uniformly bounded overlap and
uniformly bounded extension/approximation constants after scaling
to unit diameter, and their sampling sets are uniformly norming.
More precisely, for a constant $\Lambda_m$ independent of $h$ and $K$,
\[
\|p\|_{L^\infty(S(K))}
\leq\Lambda_m\max_{\boldsymbol{x}\in X(K)}|p(\boldsymbol{x})|,
\quad\forall p\in\mathbb P_m(S(K)).
\]
These assumptions imply unisolvency and the uniform reconstruction
estimate below \cite{li2024RDA}; the cardinality rule for the patches
alone is not an unisolvency assumption. For
$2\leq s\leq m+1$ and $0\leq k\leq s$, we have 
\begin{equation}\label{eqrecon}
\|z-\mathcal{R}^m(I_hz)\|_{H^k(K)}\leq C
h^{s-k}\|z\|_{H^{s}(S(K))},\quad\forall K\in\mathcal{T}_h,z\in 
H^{s}(\Omega). 
\end{equation}
In particular, the optimal rate $h^{m+1-k}$ is recovered when $z\in H^{m+1}(\Omega)$.
From now on, we abbreviate $\mathcal{R}^m(I_hz)$ to $\mathcal{R}^mz$ 
for $z\in C(\bar{\Omega})$.

\subsection{The Reconstructed Discontinuous Approximation Method}
\label{sec2.3}
Here we present the numerical method for the 
biharmonic problem based on the interior penalty discontinuous Galerkin method and the reconstructed space $\Phi_h^m$. 

We first introduce some notation commonly used in the DG 
framework. Let $e\in\mathcal{E}_h^I$ be any interior edge 
shared by $K^+,K^-\in\mathcal{T}_h$ with the unit 
outward normal vectors $\mathbf{n}^+,\mathbf{n}
^-$ along $e$, respectively. 
For any piecewise smooth scalar-valued function $v$ and 
vector-valued function $\boldsymbol
{\kappa}$, the jump operator $\jump{\,\cdot\,}$ and the average 
operator $\{\,\cdot\,\}$ are defined as
$$\jump{v}|_{e}:=v^{+}|_{e}\mathbf{n}^{+}+v^{-}|_{e}\mathbf
{n}^{-},\quad\{v\}|_{e}:=\frac{1}{2}(v^{+}|_{e}+v^{-}|_{e}),$$
$$\jump{\boldsymbol{\kappa}}|_{e}:=\boldsymbol{\kappa}^{+}|_{e}\cdot
\mathbf{n}^{+}+\boldsymbol{\kappa}^{-}
|_{e}\cdot\mathbf{n}^{-},\quad\{\boldsymbol{\kappa}\}|_{e}:=
\frac{1}{2}
(\boldsymbol{\kappa}^{+}|_{e}+\boldsymbol{\kappa}^{-}|_{e}),$$
where $v^\pm:=v|_{K^\pm},\boldsymbol{\kappa}^\pm:=\boldsymbol
{\kappa}|_{K^\pm}.$ For any boundary edge  $e\in\mathcal{E}
_h^B$, with the unit outward normal $\mathbf{n}$ along 
$e,$ the above two operators are modified as
$$\jump{v}|_e:=v|_e\mathbf{n},\ \{v\}|_e:=v|_e,\quad 
\jump{\boldsymbol{\kappa}}
| _e:=\boldsymbol{\kappa }|_e\cdot \mathbf{n},\ \{
\boldsymbol{\kappa }\} | _e: = \boldsymbol{\kappa }| _e.$$ 

Consider the following biharmonic equation
\begin{equation}\label{eq8}
	\left\{
	\begin{aligned}
		&\Delta^2\psi=F\quad\mathrm{in}\ \Omega,\\
		&\psi=g_1,\ \pn\psi=g_2\quad\mathrm{on}\ \Gamma.
	\end{aligned}
	\right.
\end{equation}
Here $F\in L^{2}(\Omega)$, and we assume that the boundary pair is
compatible:
\[
(g_1,g_2)\in\gamma(H^2(\Omega))
:=\{(\gamma_0v,\gamma_1v)\mid v\in H^2(\Omega)\}.
\]
In particular, $g_1$ and $g_2$ satisfy the polygonal corner
compatibility conditions.
We seek a numerical solution $\psi_h^m\in \Phi_h^m$ 
such that
\begin{equation}\label{eq9}
	a_{h}(\psi_h^m,\phi_h^m)=l_{h}(\phi_h^m),\quad\forall 
	\phi_h^m\in \Phi_h^m,
\end{equation}
where $a_{h}:\Phi_h^m\times \Phi_h^m\rightarrow\mathbb{R}$ and
$l_{h}:\Phi_h^m\rightarrow\mathbb{R}$ are defined by
\begin{equation*}
	\begin{aligned}
		a_h(\psi_h^m, \phi_h^m) := 
		&\sum_{K \in \mathcal{T}_h} \int_K \Delta \psi_h^m \Delta \phi_h^m
		\mathrm{d}{x} - \sum_{e \in \mathcal{E}_h} \int_e \left(  \jump{
		\nabla \psi_h^m}\aver{\Delta \phi_h^m} +  \jump{\nabla \phi_h^m} 
		\aver{\Delta 	\psi_h^m}\right)
		\mathrm{d}{s}  \\
		+ &\sum_{e \in \mathcal{E}_h} \int_e \left(\jump{\psi_h^m} \cdot
		\aver{\nabla\Delta \phi_h^m} + \jump{\phi_h^m}\cdot \aver{\nabla
		\Delta \psi_h^m} \right) \mathrm{d}{s}\\
		+ &\sum_{e \in \mathcal{E}_h} \int_e \left(\mu h_e^{-3} \jump{\psi_h^m} \cdot
		\jump{\phi_h^m}+ \mu h_e^{-1} \jump{\nabla \psi_h^m} \jump{
			\nabla \phi_h^m} \right) \mathrm{d}{s},
\\
		l_{h}(\phi_h^m):= 
		& \sum_{K\in\mathcal{T}_h}
		\int_KF\phi_h^m\mathrm{d}x+\sum_{e\in\mathcal
		{E}_h^B}\int_e\left(g_1
		\pn \Delta \phi_h^m -g_2
		\Delta \phi_h^m\right)\mathrm{d}s\quad\quad\quad\ \;\\
		+ &\sum_{e\in\mathcal{E}_h^B}\int_e \left(\mu h_e^{-3} g_1
		\phi_h^m + \mu h_e^{-1} g_2 
		\pn \phi_h^m \right) \mathrm{d}{s}.
	\end{aligned}
\end{equation*}
Here, $\mu>0$ is the penalty parameter.  The same
formulas are used as extended forms whenever the displayed element
and edge traces exist; in particular this applies to arguments in
$H^4(\mathcal T_h)$ and to the piecewise-polynomial
Hsieh--Clough--Tocher functions used below.
Denote 
\[H^s(\mathcal{T}_h):=\{\phi_h\in L^2(\Omega)\mid 
\phi_h|_{K}\in H^s(K),\forall K\in\mathcal{T}_h\}.\] 
We define the following DG-norms 
$$\|\phi_h\|_{\mathrm{DG}}^2:=\sum_{K\in\mathcal{T}_h}^{}\|
\Delta \phi_h\|_{L^2(K)}^2+\sum_{e\in\mathcal{E}_h} (h_e^
{-3}\|\jump{\phi_h}
\|_{L^2(e)}^2+h_e^{-1}\|\jump{\nabla \phi_h}\|_{L^2(e)}^2),
\quad\forall \phi_h\in H^2(\mathcal{T}_h),$$
$$\interleave \phi_h\interleave_{\mathrm{DG}}^2:=\|\phi_h\|
_{\mathrm{DG}}^2
+\sum_{e\in\mathcal{E}_h}({h_e}\|\{\Delta \phi_h\}
\|_{L^2(e)}^2
+{h_e^3}\|\{\nabla\Delta \phi_h\}
\|_{L^2(e)}^2),\quad\forall \phi_h\in H^4(\mathcal{T}_h).$$
By the inverse estimate 
and the trace inequality, we know that 
both norms are equivalent when restricted to the finite dimensional 
space $\Phi_h^m$, i.e.,  
$$\|\phi_h^m\|_{\mathrm{DG}}\leq\interleave \phi_h^m\interleave_
{\mathrm{DG}}\leq C\|\phi_h^m\|_{\mathrm{DG}},\quad\forall 
\phi_h^m\in \Phi_h^m.$$ 

The energy norm $\|\cdot\|_{\mathrm{DG}}$ is well defined on the broken space $H^2(\mathcal{T}_h)$, whereas the stronger norm $\interleave\cdot\interleave_{\mathrm{DG}}$ requires the displayed higher edge traces; elementwise $H^t$, $t>7/2$, is sufficient. We use $H^4$ regularity whenever this norm is applied to an exact function. On the discrete space $\Phi_h^m$ the two norms are equivalent, so all stability estimates below may be stated indifferently in either norm. The relationship between the energy and $L^2$ norms is given by the following lemma.
\begin{lem}\label{lem1}
There exists a constant $C>0$ such that 
\begin{equation*}
	\|\phi_h\|_{L^2(\Omega)}\leq C\|\phi_h\|_{\mathrm{DG}},\quad\forall 
	\phi_h\in H^2(\mathcal{T}_h).
\end{equation*} 
\end{lem}

\begin{proof}
Consider the second-order elliptic problem 
\[-\Delta z=\phi_h\quad \mathrm{in}\ \Omega,\quad
z=0\quad\mathrm{on}\ \Gamma,\]
where $\phi_h\in H^2(\mathcal{T}_h).$ Then $\|z\|_{H^2(\Omega)}\leq C
\|\phi_h\|_{L^2(\Omega)}$ and 
\begin{equation*}
\begin{aligned}
-(\Delta z,\phi_h)=
&-(z,\Delta \phi_h)
-\sum_{e\in\mathcal{E}_h}\int_e\nabla z\cdot \jump{\phi_h}\mathrm{d}s
+\sum_{e\in\mathcal{E}_h}\int_ez\jump{\nabla\phi_h}\mathrm{d}s\\
\leq & C\bigg(\|z\|_{L^2(\Omega)}^2
+\sum_{e\in\mathcal{E}_h}\|h_e^{\frac32}\pn z\|_{L^2(e)}^2
+\sum_{e\in\mathcal{E}_h}\|h_e^{\frac12}z\|_{L^2(e)}^2\bigg)^{\frac12}
\|\phi_h\|_{\mathrm{DG}}\\
\leq& C\|z\|_{H^2(\Omega)}\|\phi_h\|_{\mathrm{DG}}\leq C\|\phi_h\|_{L^2(\Omega)}\|\phi_h\|_{\mathrm{DG}}.
\end{aligned}
\end{equation*}
Therefore, $\|\phi_h\|_{L^2(\Omega)}\leq C\|\phi_h\|_{\mathrm{DG}}.$
Here $C$ depends on $\Omega$ and on the shape regularity of
$\mathcal{T}_h$.
\end{proof}

\begin{rmk}\label{rmk:consistency}
If the exact solution $\psi$ of $\eqref{eq8}$ belongs to $H^4(\Omega)$, integration by parts on each element yields the consistency identity
\begin{equation}\label{eqcon}
	a_{h}(\psi,\phi_h^m)=l_{h}(\phi_h^m),\quad\forall \phi_h^m
	\in \Phi_h^m,
\end{equation}
and hence the Galerkin orthogonality $a_h(\psi-\psi_h^m,\phi_h^m)=0$. 
However, as stated in Theorem \ref{thm:regularity}, the state and adjoint state of the optimal control problem may fail to belong to $H^4(\Omega)$ even when $\Omega$ is a convex polygon. Under the stated low-regularity assumptions the edge traces needed for this identity are not all available. For example, $\psi\in H^{3-\delta}(\Omega)$ with $0\leq\delta<1/2$ has an $L^2$ trace of $\Delta\psi$, but generally not of $\nabla\Delta\psi$. Thus neither \eqref{eqcon} nor state Galerkin orthogonality is used in the low-regularity analysis. We therefore adopt the medius analysis of \cite{gudi2010}, which replaces Galerkin orthogonality by residual estimates from a posteriori error analysis and compares the discontinuous space with a conforming subspace of $H_0^2(\Omega)$ through an enriching operator.
\end{rmk}

Let $V_c\subseteq H_0^2(\Omega)$ denote the Hsieh--Clough--Tocher (HCT) finite element space associated with $\mathcal{T}_h$ (cf. \cite{brenner2008book,gudi2010}). Define the enriching operator $E_h:\Phi_h^m\to V_c$ by averaging the degrees of freedom of $V_c$: if $N$ is any global degree of freedom of $V_c$ (point value or first derivative at a vertex, or normal derivative at an edge midpoint), set
\[
N(E_h\phi_h^m):=\frac{1}{|\mathcal{T}_N|}\sum_{K\in\mathcal{T}_N}N(\phi_h^m|_K)
\]
for interior degrees of freedom, where $\mathcal{T}_N$ is the patch of elements sharing $N$, and $N(E_h\phi_h^m):=0$ on $\Gamma$.

\begin{lem}[{\cite[Section~4, Eq. (4.24)]{gudi2010}}]\label{lem:enrich}
There exists a constant $C>0$ independent of $h$ such that, for all $\phi_h^m\in\Phi_h^m$,
\begin{equation}\label{eqenrich}
\sum_{K\in\mathcal{T}_h}\Bigl(h_K^{-4}\|\phi_h^m-E_h\phi_h^m\|_{L^2(K)}^2
+h_K^{-2}|\phi_h^m-E_h\phi_h^m|_{H^1(K)}^2\Bigr)
+\|E_h\phi_h^m\|_{H^2(\Omega)}^2
\leq C\|\phi_h^m\|_{\mathrm{DG}}^2.
\end{equation}
\end{lem}

Let $\widehat{\mathcal T}_h$ be the HCT subtriangulation.
Higher order derivatives of an HCT function are taken separately in
these subtriangles. Its trace on an original edge is the
one-sided trace of the polynomial adjacent to that edge.
For functions with these traces, we define
\[
T_h(d)^2:=\sum_{e\in\mathcal E_h}
\left(h_e\|\aver{\Delta d}\|_{L^2(e)}^2+
h_e^3\|\aver{\nabla\Delta d}\|_{L^2(e)}^2\right),
\]
and the following one-sided trace quantity:
\begin{equation}\label{eq:weighted-test-norm}
\begin{aligned}
D_h(d)^2:={}&
\sum_{K\in\mathcal T_h}\left(
h_K^{-4}\|d\|_{L^2(K)}^2+
h_K^{-2}|d|_{H^1(K)}^2+|d|_{H^2(K)}^2\right)\\
+&\sum_{K\in\mathcal T_h}\sum_{e\subseteq\partial K}
\left(h_e^{-3}\|d|_K\|_{L^2(e)}^2+
h_e^{-1}\|\nabla d|_K\|_{L^2(e)}^2\right.\left.+
h_e\|\Delta d|_K\|_{L^2(e)}^2+
h_e^3\|\nabla\Delta d|_K\|_{L^2(e)}^2\right).
\end{aligned}
\end{equation}
For piecewise polynomials of fixed degree on the HCT refinement,
inverse inequalities give
\begin{equation}\label{eq:macro-inverse}
D_h(d)^2\leq C
\sum_Kh_K^{-4}\|d\|_{L^2(K)}^2.
\end{equation}
Here and below the functions considered in $D_h$ are $H^2$ on
each original triangle. In particular, Lemma~\ref{lem:enrich}
and \eqref{eq:macro-inverse} imply
\begin{equation}\label{eq:enrich-test-bound}
D_h(\phi_h^m-E_h\phi_h^m)+\|E_h\phi_h^m\|_{H^2(\Omega)}
\leq C\|\phi_h^m\|_{\mathrm{DG}},\quad 
\forall\phi_h^m\in\Phi_h^m.
\end{equation}

\begin{lem}\label{lem:dual-splitting}
Let $m\geq2$, $\epsilon=\min\{m-1,2\}$, and
$z\in H^4(\Omega)\cap H_0^2(\Omega)$.
Let $I_cz\in V_c$ be its cubic HCT interpolant with clamped boundary conditions.
Then, 
\begin{equation}\label{eq:dual-splitting}
\begin{aligned}
\|z-I_cz\|_{H^2(\Omega)}+T_h(z-I_cz)&\leq Ch^2\|z\|_{H^4(\Omega)},\\
D_h(I_cz-\mathcal R^mz)&\leq Ch^\epsilon\|z\|_{H^4(\Omega)}.
\end{aligned}
\end{equation}
\end{lem}
\begin{proof}
Since $z\in H^4(\Omega)\cap H_0^2(\Omega)$, the value and gradient of $z$ are continuous, so all HCT degrees
of freedom are well defined. Its clamped boundary conditions make the boundary
degrees of freedom zero. Using polynomial reproduction and scaling on
the reference HCT macroelement, we can obtain
\[
|z-I_cz|_{H^j(\widehat{\mathcal T}_h|_K)}
\leq Ch_K^{4-j}|z|_{H^4(K)},\quad 0\leq j\leq4.
\]
For $j\leq2$ this is also the usual seminorm on $K$; for $j=4$
the fourth derivatives of the cubic interpolant vanish.
Applying the scaled trace inequality to the second and third
derivatives on boundary subtriangles gives
\[
\begin{aligned}
h_e\|\Delta (z-I_cz)\|_{L^2(e)}^2
&\leq C\Big(|z-I_cz|_{H^2(\widehat K)}^2+
h_K^2|z-I_cz|_{H^3(\widehat K)}^2\Big),\\
h_e^3\|\nabla\Delta (z-I_cz)\|_{L^2(e)}^2
&\leq C\Big(h_K^2|z-I_cz|_{H^3(\widehat K)}^2+
h_K^4|z-I_cz|_{H^4(\widehat K)}^2\Big),
\end{aligned}
\]
where $\widehat K$ is the subtriangle adjacent to $e$.
A summation over elements and edges proves the first estimate.
The $L^2$ interpolation and reconstruction estimates, together
with finite patch overlap, give
\[
\begin{aligned}
\Big(\sum_Kh_K^{-4}\|I_cz-\mathcal R^mz\|_{L^2(K)}^2\Big)^{\frac12}
\leq C(h^2\|z\|_{H^4(\Omega)}
+h^{m-1}\|z\|_{H^4(\Omega)})\leq Ch^\epsilon\|z\|_{H^4(\Omega)}.
\end{aligned}
\]
Now applying \eqref{eq:macro-inverse} finishes the proof. 
\end{proof}

\begin{lem}[\cite{georgoulis2009}]\label{lem2}
	Let $a_{h}$ be 
	defined with sufficiently large $\mu$. Then there exist constants 
	$C,c>0$ such that 
	$$|a_{h}(\psi_h^m,\phi_h^m)|\leq C\interleave \psi_h^m\interleave_
	{\mathrm{DG}} \interleave \phi_h^m\interleave_{\mathrm{DG}},
	\quad\forall \psi_h^m,\phi_h^m\in \Phi_h^m,$$
	$$a_{h}(\phi_h^m,\phi_h^m)\geq c \interleave \phi_h^m
	\interleave_{\mathrm{DG}}^2, \quad\forall \phi_h^m\in \Phi_h^m.$$
	In particular, the same bounds hold with $\|\cdot\|_{\mathrm{DG}}$ in place of $\interleave\cdot\interleave_{\mathrm{DG}}$.
\end{lem}

For later use, we record two residual identities that require
only the weak solution $\psi\in H^2(\Omega)$ of \eqref{eq8}.
Write
\[
\mathcal F_{\phi_h^m}(z):=l_h(z)-a_h({\phi_h^m},z),
\quad \forall {\phi_h^m}\in\Phi_h^m.
\]
If $z\in H_0^2(\Omega)$ has the higher original-edge traces
appearing in $T_h(z)$, the variational formulation for the biharmonic equation and the prescribed
boundary traces give
\begin{equation}\label{eq:conforming-residual}
\mathcal F_{\phi_h^m}(z)
=(\Delta (\psi-{\phi_h^m}),\Delta z)
-\sum_{e\in\mathcal E_h}(\jump{\nabla (\psi-{\phi_h^m})},\aver{\Delta z})_e
+\sum_{e\in\mathcal E_h}(\jump{\psi-{\phi_h^m}},\aver{\nabla\Delta z})_e.
\end{equation}
The edge inner products include a dot product when the arguments
are vectors. We note that no higher trace of $\psi$ occurs in this identity.
It follows that
\begin{equation}\label{eq:conforming-residual-bound}
|\mathcal F_{\phi_h^m}(z)|
\leq C\|\psi-{\phi_h^m}\|_{\mathrm{DG}}
\bigl(\|z\|_{H^2(\Omega)}+T_h(z)\bigr),\quad\forall z\in H_0^2(\Omega).
\end{equation}
In particular, inverse traces for $z\in V_c$ imply
\begin{equation}\label{eqN2}
|\mathcal F_{\phi_h^m}(z)|\leq
C\|\psi-{\phi_h^m}\|_{\mathrm{DG}}\|z\|_{H^2(\Omega)},
\quad\forall z\in V_c.
\end{equation}
We note that the terms containing $g_1,g_2$ in $l_h(z)$ are needed here:
they turn the boundary terms into traces of $\psi-{\phi_h^m}$.

\begin{lem}\label{lem:efficiency}
Let $\psi\in H^2(\Omega)$ be the weak solution of \eqref{eq8},
$F\in L^2(\Omega)$, and ${\phi_h^m}\in\Phi_h^m$. Define
\[
\operatorname{Osc}(F)^2:=
\sum_K h_K^4
\inf_{p\in\mathbb P_{\max\{m-2,0\}}(K)}
\|F-p\|_{L^2(K)}^2.
\]
Then
\begin{equation}\label{eq:residual-efficiency}
\sum_Kh_K^4\|F-\Delta^2{\phi_h^m}\|_{L^2(K)}^2+
\sum_{e\in\mathcal E_h^I}
\Big(h_e\|\jump{\Delta {\phi_h^m}}\|_{L^2(e)}^2+
h_e^3\|\jump{\nabla\Delta {\phi_h^m}}\|_{L^2(e)}^2\Big)\leq C\bigl(\|\psi-{\phi_h^m}\|_{\mathrm{DG}}
+\operatorname{Osc}(F)\bigr)^2.
\end{equation}
The boundary mismatches satisfy
\[
\sum_{e\in\mathcal E_h^B}
\Big(h_e^{-3}\|g_1-{\phi_h^m}\|_{L^2(e)}^2+
h_e^{-1}\|g_2-\pn {\phi_h^m}\|_{L^2(e)}^2\Big)
\leq\|\psi-{\phi_h^m}\|_{\mathrm{DG}}^2.
\]
\end{lem}
\begin{proof}
The interior estimates are those of
\cite[Lemma~4.2]{gudi2010}. Their test functions belong to
$H_0^2$ of an element or a two-element patch, extended by zero
to $\Omega$, so the proof is unchanged for  compatible
inhomogeneous boundary data here. We give the local estimates
to make this point explicit.
Let $F_K$ be the elementwise $L^2$ projection onto
$\mathbb P_{\max\{m-2,0\}}(K)$ and let
$r_K=F_K-\Delta^2{\phi_h^m}$. Testing with $b_K^2r_K$, where $b_K$ is
the element bubble, and using the weak formulation for biharmonic equations gives
\[
(r_K,b_K^2r_K)_K
=(\Delta(\psi-{\phi_h^m}),\Delta(b_K^2r_K))_K
+(F_K-F,b_K^2r_K)_K.
\]
We can conclude from the norm equivalence for polynomials and the scaling argument that
\[
h_K^2\|F-\Delta^2{\phi_h^m}\|_{L^2(K)}
\leq C\left(\|\Delta(\psi-{\phi_h^m})\|_{L^2(K)}
+h_K^2\|F-F_K\|_{L^2(K)}\right).
\]
For an interior edge $e$, let $\omega_e$ be its two-element
patch. The clamped edge-bubble liftings of
$\jump{\Delta {\phi_h^m}}$ and $\jump{\nabla\Delta {\phi_h^m}}$ give, respectively,
\[
\begin{aligned}
h_e^{\frac12}\|\jump{\Delta {\phi_h^m}}\|_{L^2(e)}
&\leq C\Big(\|\Delta(\psi-{\phi_h^m})\|_{L^2(\omega_e)}
+\sum_{K\subseteq\omega_e}h_K^2
\|F-\Delta^2{\phi_h^m}\|_{L^2(K)}\Big),\\
h_e^{\frac32}\|\jump{\nabla\Delta {\phi_h^m}}\|_{L^2(e)}
&\leq C\Big(\|\Delta(\psi-{\phi_h^m})\|_{L^2(\omega_e)}
+\sum_{K\subseteq\omega_e}h_K^2
\|F-\Delta^2{\phi_h^m}\|_{L^2(K)}
+h_e^{\frac12}\|\jump{\Delta {\phi_h^m}}\|_{L^2(e)}\Big).
\end{aligned}
\]
These follow by integrating the polynomial ${\phi_h^m}$ twice and using
the weak formulation with the liftings; their scaled bounds are
precisely the edge-bubble estimates in the cited lemma.
Squaring and summing, with bounded patch overlap, proves
\eqref{eq:residual-efficiency}. The boundary estimate follows
directly from the definition of the DG norm and the prescribed
traces of $\psi$; it does not require the reconstruction of boundary bubbles.
\end{proof}

If $d$ is $H^2$ in each original element and is a piecewise polynomial
in its HCT refinement, integration by parts of ${\phi_h^m}$ gives
\begin{equation}\label{eq:polynomial-residual}
\begin{aligned}
\mathcal F_{\phi_h^m}(d)={}&(F-\Delta^2{\phi_h^m},d)
-\sum_{e\in\mathcal E_h^I}(\jump{\Delta {\phi_h^m}},\aver{\nabla d})_e
+\sum_{e\in\mathcal E_h^I}(\jump{\nabla\Delta {\phi_h^m}},\aver d)_e\\
-&\sum_{e\in\mathcal E_h}(\jump{\nabla (\psi-{\phi_h^m})},\aver{\Delta d})_e
+\sum_{e\in\mathcal E_h}(\jump{\psi-{\phi_h^m}},\aver{\nabla\Delta d})_e\\
+&\sum_{e\in\mathcal E_h}\mu
\left(h_e^{-3}(\jump{\psi-{\phi_h^m}},\jump d)_e+
h_e^{-1}(\jump{\nabla (\psi-{\phi_h^m})},\jump{\nabla d})_e\right).
\end{aligned}
\end{equation}
There are no terms at the internal HCT edges, since $d$ and its
gradient are continuous there and ${\phi_h^m}$ is a single polynomial in
each original triangle. On a boundary edge, the terms containing
$g_1,g_2$ have been combined with ${\phi_h^m},\pn {\phi_h^m}$ to form $e_{\phi_h^m},\pn e_{\phi_h^m}$.
The weighted Cauchy--Schwarz inequality and Lemma~\ref{lem:efficiency} imply
\begin{equation}\label{eq:polynomial-residual-bound}
|\mathcal F_{\phi_h^m}(d)|
\leq C\bigl(\|\psi-{\phi_h^m}\|_{\mathrm{DG}}+
\operatorname{Osc}(F)\bigr)D_h(d).
\end{equation}

\begin{thm}\label{thm1}
	Let $2\leq s\leq m+1$ and let $\psi\in H^{s}(\Omega)$ be the weak solution of 
	$\eqref{eq8}$, that is, $\psi=g_1$ and $\pn\psi=g_2$ on $\Gamma$ and
	$a(\psi,\phi)=(F,\phi)$ for all $\phi\in H_0^2(\Omega)$. 
	Let $\psi_h^m \in \Phi_h^m$ be the numerical 
	solution of $\eqref{eq9}$ with $m\geq 2$ and $\mu$ 
	sufficiently large. Then there exists a constant 
	$C>0$ independent of $h$ such that
	\begin{equation}\label{eq:thm1energy}
	\|\psi - \psi_h^m\|_{\mathrm{DG}} \leq C 
	\bigl(h^{s-2} \|\psi\|_{H^{s}(\Omega)}+\mathrm{Osc}(F)\bigr).
	\end{equation}
	In particular, if $\psi\in H^{m+1}(\Omega)$ and
	$\mathrm{Osc}(F)\leq Ch^{m-1}\|\psi\|_{H^{m+1}(\Omega)}$, then
	$\|\psi - \psi_h^m\|_{\mathrm{DG}} \leq C 
	h^{m-1} \|\psi\|_{H^{m+1}(\Omega)}.$
    
	Moreover, if every interior angle of $\Omega$ is less than $126.283$ degrees, then
	\begin{equation}\label{eq:thm1L2}
	\|\psi - \psi_h^m\|_{L^2(\Omega)} \leq C 
	h^{\ep}\bigl(h^{s-2}\|\psi\|_{H^{s}(\Omega)}
	+\mathrm{Osc}(F)\bigr),
	\end{equation}
	where $\ep=\min\{m-1,2\}$. Under the preceding oscillation
	assumption, this gives
	$\|\psi - \psi_h^m\|_{L^2(\Omega)} \leq C 
	h^{m+\ep-1} \|\psi\|_{H^{m+1}(\Omega)}.$
\end{thm}
\begin{proof}
Let $\phi_h^m\in\Phi_h^m$ and $\varphi=\psi_h^m-\phi_h^m$. The coercivity of 
$a_h(\cdot,\cdot)$ and
the discrete equation give
\[
c\|\varphi\|_{\mathrm{DG}}^2\leq a_h(\psi_h^m-\phi_h^m,\varphi)=l_h(\varphi)-a_h({\phi_h^m},\varphi)=
\mathcal F_{\phi_h^m}(\varphi)
=\mathcal F_{\phi_h^m}(E_h\varphi)
+\mathcal F_{\phi_h^m}(\varphi-E_h\varphi).
\]
Use \eqref{eqN2}, \eqref{eq:enrich-test-bound} and
\eqref{eq:polynomial-residual-bound} to obtain
\[
\|\varphi\|_{\mathrm{DG}}
\leq C\bigl(\|\psi-\phi_h^m\|_{\mathrm{DG}}
+\operatorname{Osc}(F)\bigr).
\]
The triangle inequality yields the medius estimate
\begin{equation}\label{eq:medius-best}
\|\psi-\psi_h^m\|_{\mathrm{DG}}
\leq C\Big(\inf_{\phi_h^m\in\Phi_h^m}\|\psi-\phi_h^m\|_{\mathrm{DG}}
+\operatorname{Osc}(F)\Big).
\end{equation}
The reconstruction estimate \eqref{eqrecon}, scaled traces for
the value and first derivative, and bounded patch overlap give
$\|\psi-\mathcal R^m\psi\|_{\mathrm{DG}}
\leq Ch^{s-2}\|\psi\|_{H^s(\Omega)}$. This proves
\eqref{eq:thm1energy}. Notice that
$\operatorname{Osc}(F)\leq Ch^2\|F\|_{L^2(\Omega)}$
for arbitrary $L^2$ data.

For the $L^2$ norm error estimate, we set $e=\psi-\psi_h^m$ and solve the following equation
\[
\Delta^2z=e\quad\mathrm{in}\ \Omega,\quad z\in H_0^2(\Omega).
\]
The angle restriction implies
$\|z\|_{H^4(\Omega)}\leq C\|e\|_{L^2(\Omega)}$.
Since $\psi\in H^2(\Omega)$ and $z\in H^4(\Omega)\cap H_0^2(\Omega)$, Green's
formula and the weak formulation imply
\[
(\psi,\Delta^2z)
=(F,z)+(g_1,\pn\Delta z)_\Gamma-(g_2,\Delta z)_\Gamma
=l_h(z),
\]
where we used $z=\pn z=0$ on $\Gamma$.
Integration and summation over elements give 
$a_h(\psi_h^m,z)=(\psi_h^m,\Delta^2z)$. Consequently,
using only the discrete equation with $\mathcal R^mz$, we obtain
\begin{equation*}
\|e\|_{L^2(\Omega)}^2
=\mathcal F_{\psi_h^m}(z)
=\mathcal F_{\psi_h^m}(z-\mathcal R^mz)=
\mathcal F_{\psi_h^m}(z-I_cz)+\mathcal F_{\psi_h^m}(I_cz-\mathcal R^mz).
\end{equation*}

Apply \eqref{eq:conforming-residual-bound} to  
 $z-I_cz$ with clamped boundary conditions; and apply equation
\eqref{eq:polynomial-residual-bound} to the piecewise
polynomial $I_cz-\mathcal R^mz$. Thus, Lemma \ref{lem:dual-splitting} implies
\[
\begin{aligned}
\|e\|_{L^2(\Omega)}^2
&\leq C\|e\|_{\mathrm{DG}}
\bigl(\|z-I_cz\|_{H^2(\Omega)}+T_h(z-I_cz)\bigr)
+C\bigl(\|e\|_{\mathrm{DG}}+\operatorname{Osc}(F)\bigr)D_h(I_cz-\mathcal R^mz)\\
&\leq Ch^\epsilon
\bigl(\|e\|_{\mathrm{DG}}+\operatorname{Osc}(F)\bigr)
\|z\|_{H^4(\Omega)}.
\end{aligned}
\]
Divide by $\|e\|_{L^2(\Omega)}$ when it is nonzero and apply
\eqref{eq:thm1energy} to obtain \eqref{eq:thm1L2}.
The simplified estimates follow immediately from the assumed
oscillation bound.
\end{proof}

For a sufficiently smooth exact solution, the strong DG norm
can be estimated without an oscillation term as in Theorem \ref{thm11}.
\begin{thm}\label{thm11}
Let $m\geq3$, $4\leq t\leq m+1$, and let
$\psi\in H^t(\Omega)$ solve \eqref{eq8}. If $\psi_h^m$
solves \eqref{eq9} and the penalty is sufficiently large, then
\[
\interleave\psi-\psi_h^m\interleave_{\mathrm{DG}}
\leq Ch^{t-2}\|\psi\|_{H^t(\Omega)}.
\]
In particular, $t=m+1$ gives order $h^{m-1}$, while $t=4$
gives order $h^2$ for every $m\geq3$.
\end{thm}
\begin{proof}
The reconstruction estimate and scaled traces of its second
and third derivatives give
\[
\interleave\psi-\mathcal R^m\psi\interleave_{\mathrm{DG}}
\leq Ch^{t-2}\|\psi\|_{H^t(\Omega)}.
\]
The formula for $a_h$ is continuous in the strong DG norm
whenever its traces exist, by the weighted Cauchy--Schwarz inequality.
Since $\psi\in H^4(\Omega)$, consistency is valid. Put
$\mathcal R^m\psi-\psi_h^m\in\Phi_h^m$. Coercivity,
consistency, and discrete norm equivalence give
\[
\begin{aligned}
c\|\mathcal R^m\psi-\psi_h^m\|_{\mathrm{DG}}^2&
\leq a_h(\mathcal R^m\psi-\psi_h^m,\mathcal R^m\psi-\psi_h^m)
=a_h(\mathcal R^m\psi-\psi,\mathcal R^m\psi-\psi_h^m)\\
&\leq C\interleave\mathcal R^m\psi-\psi\interleave_{\mathrm{DG}}
\|R^m\psi-\psi_h^m\|_{\mathrm{DG}}\leq Ch^{t-2}\|\psi\|_{H^t(\Omega)}\|\mathcal R^m\psi-\psi_h^m\|_{\mathrm{DG}}.
\end{aligned}
\]
The triangle inequality and discrete norm equivalence finish
the proof.
\end{proof}

\begin{rmk}\label{rmk1}
	Here we use the inverse estimate to furnish an upper bound 
	for the piecewise norm $|\psi-\psi_h^m|_{H^k(\mathcal T_h)}$.
	Assume the angle and oscillation hypotheses of the simplified
	$L^2$ estimate in Theorem \ref{thm1}. For $0\leq k<m+\ep-1$, consider another 
	reconstruction operator $R^m:L^2(\Omega)\rightarrow\mathbb{P}_m(\mathcal{T}_h)$, 
	the $L^2$-projection to  
	the discontinuous piecewise $m$-th degree polynomial space. 
	We have 
	\[|\psi-R^m\psi_h^m|_{H^k(\mathcal T_h)}\leq|\psi-R^m\psi|_{H^k(\mathcal T_h)}+
	|R^m(\psi-\psi_h^m)|_{H^k(\mathcal T_h)}.\]
	In addition, we have 
	$|\psi-R^m\psi|_{H^k(\mathcal T_h)}\leq Ch^{m+1-k}\|\psi\|_{H^{m+1}(\Omega)}$ and 
	$$|R^m(\psi-\psi_h^m)|_{H^k(\mathcal T_h)}\leq Ch^{-k}\|R^m(\psi-\psi_h^m)\|_
	{L^2(\Omega)}\leq Ch^{-k}\|\psi-\psi_h^m\|_
	{L^2(\Omega)}$$
	by the inverse estimate. Thus
	$|\psi-R^m\psi_h^m|_{H^k(\mathcal T_h)}
	\leq Ch^{m+\ep-1-k}
	\|\psi\|_{H^{m+1}(\Omega)}.$
	Notice that $R^m\psi_h^m=\psi_h^m,$ we can derive the 
	desired results.
\end{rmk}

\subsection{The Discrete First Order Optimality System.} \label{sec2.5}
So far, we have established the numerical scheme for the 
biharmonic equation and its convergence analysis under low regularity. We now
give a discretization of the optimality system  
$\eqref{eq6}$. We denote $(\psi^*,
u^*,\eta^*)$ as the exact solution. 
Take $r\geq1$ and let $U_{h}^r\subseteq U_{ad}$ be a finite dimensional 
space, whose restriction on each edge $\Gamma_i$
$(i=1,2,\cdots ,l)$ is $P_r(\Gamma_{i})$, the  
continuous piecewise polynomial space of degree $r$ on the mesh of
$\Gamma_i$. In the case $U=H^1(\Gamma)$ we impose
\[
v_h^r(P)=0,\quad\forall P\in\mathcal V_\Gamma,\  v_h^r\in U_h^r,
\]
so that $U_h^r$ is a conforming subspace of the corner-compatible
space $U_{ad}$. Therefore, in the case $U=H^1(\Gamma)$ we adopt a continuous discrete control space $U_h^r$. We assume that all polygonal vertices are nodes of
the boundary mesh. In the $L^2(\Gamma)$ case continuity is imposed
within each open macro-edge only; no cross-corner continuity is
required.

When the exact solution is sufficiently regular (e.g. $\psi^*,\eta^*\in H^4(\Omega)$), the consistency identity \eqref{eqcon} implies
\begin{equation}\label{eq12}
	\left\{
	\begin{aligned}
		&a_{h}(\psi^{\ast},\phi_h^m)=l_{h}[f,u^*]
		(\phi_h^m),\quad\forall \phi_h^m\in \Phi_h^m,\\
		&a_{h}(\xi_h^m,\eta^{\ast})=l_{h}[\psi^*-\psi_d
		,0](\xi_h^m),\quad\forall \xi_h^m\in \Phi_h^m,\\
		&\alpha (u^*,v_h^r)_U=-l_h[0,v_h^r](\eta^*),\quad\forall v_h^r\in U_{h}^r.
	\end{aligned}
	\right.
\end{equation}
In general, we do not rely on \eqref{eq12}; the discrete scheme is defined directly from the optimality system \eqref{eq6} as follows.
\begin{equation}\label{eq13}
	\left\{
	\begin{aligned}
		&a_{h}(\psi_h^{m\ast},\phi_h^m)=l_{h}[ f,u_h^{r*}]
		(\phi_h^m),\quad\forall \phi_h^m\in \Phi_h^m,\\
		&a_{h}(\xi_h^m,\eta_h^{m\ast})=l_{h}[\psi_h^{m*}-\psi_d
		,0](\xi_h^m),\quad\forall \xi_h^m\in \Phi_h^m,\\
		&\alpha (u_h^{r*},v_h^r)_U=-l_h[0,v_h^r](\eta_h^{m*}),
		\quad\forall v_h^r\in U_{h}^r.
	\end{aligned}
	\right.
\end{equation}
Here for $\phi_h\in H^2(\mathcal{T}_h),$ 
\begin{equation*}
	\begin{aligned}
		l_{h}[t,v](\phi_h):= 
		\sum_{K\in\mathcal{T}_h}
		\int_Kt\phi_h\mathrm{d}x-\sum_{e\in\mathcal
			{E}_h^B}\int_e v
		\Delta \phi_h\mathrm{d}s
		+ \sum_{e\in\mathcal{E}_h^B}\int_e 
		\mu h_e^{-1} v 
		\pn \phi_h \mathrm{d}{s}.
	\end{aligned}
\end{equation*}
The functional $l_h[t,v](\phi_h)$ is linear in the data pair $(t,v)$ and linear in the test function $\phi_h$.
We use the same formula for exact functions $\phi$ whenever
$\Delta\phi|_\Gamma$ and $\pn\phi|_\Gamma$ belong to
$L^2(\Gamma)$; in particular this extension is valid for
$\phi\in H^3(\Omega)$.

\begin{rmk}\label{rmk3}
We mention here that \eqref{eq13} is actually the first order optimality system 
of a discrete optimal control problem
\begin{equation*}
	\begin{aligned}
		&\min_{u_h^r\in U_{h}^r}\frac{1}{2}\|\psi_h^m-\psi_d\|_{L^2(\Omega)}^2+\frac\alpha2
		\|u_h^r\|_U^2,&\\
		\mathrm{s.t.}\ 
		&\psi_h^m\in \Phi_h^m:\quad a_h(\psi_h^m,\phi_h^m)=l_h[f,u_h^r](\phi_h^m),\quad\forall \phi_h^m\in \Phi_h^m.
	\end{aligned}
\end{equation*}
To derive the optimality condition of this problem, 
we separate
the constraint equation into two independent parts as
\begin{flalign*}
	a_h(\psi_h^f,\phi_h^m)=l_h[f,0](\phi_h^m),\quad 
	a_h(\psi_h^0(u_h^r),\phi_h^m)=l_h[0,u_h^r](\phi_h^m).
\end{flalign*}
Denote the optimal solution to this problem by $(\psi_h^{m*},
u_h^{r*})$, then 
$\psi_h^{0*}:=\psi_h^{0}(u_h^{r*})=\psi_h^{m*}-\psi_h^{f}$.  
Define the adjoint variable $\eta_h^{m*}\in \Phi_h^m$ such that 
\[a_h(\xi_h^m,\eta_h^{m*})=l_h[\psi_h^{m*}-\psi_d,0](\xi_h^m),\quad\forall
\xi_h^m\in \Phi_h^m.\]
Denote the cost functional as $J(\psi_h^0(u_h^r),u_h^r),$ we have 
$$J'(\psi_h^{0*},u_h^{r*})(v_h^r)=(\psi_h^{0*}+\psi_h^f-\psi_d,
\psi_h^{0'}(u_h^{r*})(v_h^r))+\alpha(u_h^{r*},v_h^r)_U=0,\quad\forall v_h^r\in U_{h}^r.$$
Notice that $\psi_h^{0'}(u_h^{r*})(v_h^r)=\psi_h^0(v_h^r)$ 
for any $v_h^r\in U_{h}^r$. Therefore,  
\[a_h(\psi_h^0(v_h^r),\eta_h^{m*})+\alpha(u_h^{r*},v_h^r)_U=0,
\quad\forall v_h^r\in U_{h}^r.\] 
Notice that 
\begin{equation*}
	\begin{aligned}
		a_h(\psi_h^0(v_h^r),\eta_h^{m*})=l_h[0,v_h^r](\eta_h^{m*}),
		\quad\forall v_h^r\in U_{h}^r,
	\end{aligned}
\end{equation*}
we have 
\[\alpha(u_h^{r*},v_h^r)_U=-l_h[0,v_h^r](\eta_h^{m*}),
\quad\forall v_h^r\in U_{h}^r,\]
which is exactly the last equation in \eqref{eq13}.
\end{rmk}

\section{A Priori Error Estimates}\label{sec3}
At first, we define some projection operators on the boundary of domain and give some error estimates.
\begin{lem}\label{lem3}
Define $P_hu^*\in U_{h}^r$ by 
\[\alpha(P_hu^*,v_h^r)_U=-l_h[0,v_h^r](\eta^*),\quad\forall v_h^r\in U_{h}^r.\]
Assume $u^*\in H_{\rm pw}^{r+1}(\Gamma)$ and, in the
$H^1(\Gamma)$ case, $u^*\in U_{ad}$. Then
\begin{equation}\label{eq:Phu}
\|u^*-P_hu^*\|_{L^2(\Gamma)}\leq Ch^{r+1}\|u^*\|_{H_{\rm pw}^{r+1}(\Gamma)}.
\end{equation}
Moreover, if $U=H^1(\Gamma)$,
\begin{equation}\label{eq:PhuH1}
\|u^*-P_hu^*\|_{H^1(\Gamma)}\leq Ch^{r}\|u^*\|_{H_{\rm pw}^{r+1}(\Gamma)}.
\end{equation}
\end{lem}

\begin{proof}
Since $\eta^*\in H_0^2(\Omega)$, we have $\pn\eta^*=0$ on $\Gamma$, and therefore
$-l_h[0,v](\eta^*)=(\Delta\eta^*,v)_{\Gamma}$ for all $v\in U_{ad}$. 
Combined with the optimality condition $\alpha(u^*,v)_U=(\Delta\eta^*,v)_{\Gamma}$, this yields
\begin{equation}\label{eq:Phu-orth}
\alpha(P_hu^*,v_h^r)_U=\alpha(u^*,v_h^r)_U,\quad\forall v_h^r\in U_{h}^r.
\end{equation}
Thus, $P_hu^*$ is the orthogonal projection of $u^*$ onto $U_{h}^r$ with respect to the inner product of $U$.
Let $I_h^uu^*\in U_{h}^r$ denote an interpolant of $u^*$.

If $U=L^2(\Gamma)$, then \eqref{eq:Phu-orth} is the $L^2(\Gamma)$-orthogonal projection, and C\'ea's lemma gives
\[
\|u^*-P_hu^*\|_{L^2(\Gamma)}\leq\|u^*-I_h^uu^*\|_{L^2(\Gamma)}\leq Ch^{r+1}\|u^*\|_{H_{\rm pw}^{r+1}(\Gamma)}.
\]
If $U=H^1(\Gamma)$, write
\[
(w,v)_{H^1(\Gamma)}:=(w,v)_{\Gamma}+(\nabla_\Gamma w,\nabla_\Gamma v)_{\Gamma}.
\]
Then \eqref{eq:Phu-orth} is the Ritz projection of $u^*$ onto $U_{h}^r$. C\'ea's lemma yields \eqref{eq:PhuH1}
\[
\|u^*-P_hu^*\|_{H^1(\Gamma)}\leq C\|u^*-I_h^uu^*\|_{H^1(\Gamma)}\leq Ch^{r}\|u^*\|_{H_{\rm pw}^{r+1}(\Gamma)}.
\]
To improve the $L^2(\Gamma)$-norm estimate, we use a duality argument on $\Gamma$. 
Let $e_u:=u^*-P_hu^*$ and let $w\in U_{ad}$ solve
\begin{equation}\label{eq:dual-Gamma}
(w,v)_{H^1(\Gamma)}=(e_u,v)_{\Gamma},\quad\forall v\in U_{ad}.
\end{equation}
This is a collection of one-dimensional Dirichlet problems
$-w''+w=e_u$ on $\Gamma_i$, with zero values at their
endpoints. Consequently
$w\in H_{\rm pw}^2(\Gamma)\cap U_{ad}$ and
$\|w\|_{H_{\rm pw}^2(\Gamma)}
\leq C\|e_u\|_{L^2(\Gamma)}$.
The edgewise nodal interpolant $I_h^u\phi$ belongs to $U_h^r$ because
it preserves the zero vertex values, and
\begin{equation}\label{eq:interp-Gamma}
\|v-I_h^uv\|_{H^1(\Gamma)}
\leq Ch\|v\|_{H_{\rm pw}^2(\Gamma)},\quad\forall v\in U_{ad}.
\end{equation}
Taking $v=e_u$ in \eqref{eq:dual-Gamma} and using the Galerkin orthogonality we obtain
\begin{equation*}
\begin{aligned}
\|e_u\|_{L^2(\Gamma)}^2
&=(e_u,w)_{H^1(\Gamma)}
=(e_u,w-I_h^uw)_{H^1(\Gamma)}
\leq\|e_u\|_{H^1(\Gamma)}\|w-I_h^uw\|_{H^1(\Gamma)}\\
&\leq C\|e_u\|_{H^1(\Gamma)}\,h\|w\|_{H_{\rm pw}^2(\Gamma)}
\leq Ch\|e_u\|_{H^1(\Gamma)}\|e_u\|_{L^2(\Gamma)}.
\end{aligned}
\end{equation*}
Hence $\|e_u\|_{L^2(\Gamma)}\leq Ch\|e_u\|_{H^1(\Gamma)}$, which together with \eqref{eq:PhuH1} proves \eqref{eq:Phu}.
\end{proof}

\begin{rmk}\label{cor:control-projection}
The projection $P_h$ also satisfies the fractional-order estimates
\[
\|u^*-P_hu^*\|_{L^2(\Gamma)}
\leq Ch^t\|u^*\|_{H_{\rm pw}^t(\Gamma)},
\quad 0\leq t\leq r+1,
\]
when $U=L^2(\Gamma)$. When $U=H^1(\Gamma)$ and
$u^*\in H_{\rm pw}^t(\Gamma)\cap U_{ad}$, $1\leq t\leq r+1$, we have 
\[
\|u^*-P_hu^*\|_{H^1(\Gamma)}
\leq Ch^{t-1}\|u^*\|_{H_{\rm pw}^t(\Gamma)},\quad
\|u^*-P_hu^*\|_{L^2(\Gamma)}
\leq Ch^t\|u^*\|_{H_{\rm pw}^t(\Gamma)}.
\]
Notice that 
the first two estimates are the standard edgewise approximation
estimates for the $L^2$ and Ritz projections, obtained by real
interpolation from the integer-order bounds. The last estimate
follows from the edgewise Dirichlet duality argument used in
Lemma \ref{lem3}.
\end{rmk}

We write $(\psi,u,\eta)$ and $(\psi_h^m,u_h^r,\eta_h^m)$ for the
continuous and discrete optimal solutions, i.e. \eqref{eq6}
and \eqref{eq13}. Set
$e_\psi=\psi-\psi_h^m$, $e_u=u-u_h^r$, and
$e_\eta=\eta-\eta_h^m$. All constants may depend on the fixed
domain, $\alpha$, the polynomial degrees, the penalty parameter,
and the uniform mesh and reconstruction constants, but not on $h$.

For discrete $q_h$, inverse traces give
\begin{equation}\label{eq:boundary-functional}
|l_h[0,v](\phi_h^m)|\leq Ch^{-\frac12}\|v\|_{L^2(\Gamma)}
\|\phi_h^m\|_{\mathrm{DG}},
\quad \forall v\in L^2(\Gamma), \phi_h^m\in\Phi_h^m.
\end{equation}
For any $\phi_h\in H^4(\mathcal{T}_h)$ whose strong DG norm is defined, the boundary
weights in that norm directly give
\begin{equation}\label{eq:boundary-functional-strong}
|l_h[0,v](\phi_h)|\leq Ch^{-\frac12}\|v\|_{L^2(\Gamma)}
\interleave \phi_h\interleave_{\mathrm{DG}},
\quad \forall v\in L^2(\Gamma), \phi_h\in H^4(\mathcal{T}_h).
\end{equation}

Define the two discrete elliptic projections by
\begin{equation}\label{eq:elliptic-projections}
\begin{aligned}
a_h(P_h\psi,\phi_h^m)=l_h[f,u](\phi_h^m),\quad
a_h(\phi_h^m,P_h\eta)=(\psi-\psi_d,\phi_h^m),
\quad \forall\phi_h^m\in\Phi_h^m.
\end{aligned}
\end{equation}
These definitions require only $f,\psi-\psi_d\in L^2(\Omega)$
and $u\in L^2(\Gamma)$.

\subsection{Estimates under high regularity}
In this subsection, we provide the error estimate for optimal control problems with higher regularity.
\begin{thm}\label{thm2}
Let the continuous and discrete solutions satisfy \eqref{eq6}
and \eqref{eq13}. Take $r=m-1$, $m\geq2$, and assume
$\psi,\eta\in H^{m+1}(\Omega)$ and
$u\in H_{\rm pw}^{m}(\Gamma)\cap U_{ad}$. Suppose also that
\[
\operatorname{Osc}(f)+\operatorname{Osc}(\psi-\psi_d)
\leq C_{\rm osc}h^{m-1}.
\]
Then, for a constant depending on these solution norms and
$C_{\rm osc}$, we have 
\[
\|u-u_h^r\|_U\leq Ch^{m-\frac32},\quad
\|\psi-\psi_h^m\|_{\mathrm{DG}}\leq Ch^{m-2},\quad
\|\eta-\eta_h^m\|_{\mathrm{DG}}\leq Ch^{m-\frac32}.
\]
\end{thm}
\begin{proof}
We first recall the following energy norm error estimates stated in Theorem~\ref{thm1}
\begin{equation}\label{eq:projection-energy}
\|\psi-P_h\psi\|_{\mathrm{DG}}
+\|\eta-P_h\eta\|_{\mathrm{DG}}\leq C h^{m-1}.
\end{equation}
In particular, Lemma~\ref{lem1} implies
$\|\psi-P_h\psi\|_{L^2(\Omega)}\leq Ch^{m-1}$.
We note that no angle restriction is used in this step.

Subtracting the discrete equations from
\eqref{eq:elliptic-projections} gives
\begin{equation}\label{state1}
\begin{aligned}
a_h(\psi_h^m-P_h\psi,\phi_h^m)=l_h[0,u_h^r-u](\phi_h^m),\quad
a_h(\phi_h^m,\eta_h^m-P_h\eta)=(\psi_h^m-\psi,\phi_h^m).
\end{aligned}
\end{equation}
The coercivity of $a_h$, \eqref{eq:boundary-functional}, and
Lemma~\ref{lem1} therefore give
\begin{equation}\label{eqdg}
\begin{aligned}
\|\psi_h^m-P_h\psi\|_{\mathrm{DG}}
\leq Ch^{-\frac12}\|u_h^r-u\|_{L^2(\Gamma)},\quad 
\|\eta_h^m-P_h\eta\|_{\mathrm{DG}}
\leq C\|\psi_h^m-\psi\|_{L^2(\Omega)}.
\end{aligned}
\end{equation}

We also need the estimate of the boundary functional on the adjoint projection
error, i.e., $l_h[0,v](P_h\eta-\eta)$. Split
$P_h\eta-\eta=(P_h\eta-\mathcal R^m\eta)+(\mathcal R^m\eta-\eta)$.
For the estimate of the discrete part $(P_h\eta-\mathcal R^m\eta)$ we use \eqref{eq:boundary-functional} and
\eqref{eq:projection-energy}. For the second part, the scaled
trace inequality and \eqref{eqrecon} give
\[
\begin{aligned}
\|\Delta(\eta-\mathcal R^m\eta)\|_{L^2(\Gamma)}
+\Big(\sum_{e\in\mathcal E_h^B}
h_e^{-2}\|\pn(\eta-\mathcal R^m\eta)\|_{L^2(e)}^2\Big)^{\frac12}\leq Ch^{m-\frac32}\|\eta\|_{H^{m+1}(\Omega)}.
\end{aligned}
\]
We note that this trace estimate needs only $\eta\in H^3(\Omega)$ when $m=2$.
Consequently,
\begin{equation}\label{eq:adjoint-boundary-projection}
|l_h[0,v](P_h\eta-\eta)|
\leq C h^{m-\frac32}\|v\|_{L^2(\Gamma)},
\quad\forall v\in L^2(\Gamma).
\end{equation}
The data-dependent constant includes the adjoint oscillation
bound in \eqref{eq:projection-energy}.

Let $P_hu$ be defined in  Lemma~\ref{lem3},
and define $Q_h\psi\in\Phi_h^m$ by
$$a_h(Q_h\psi,\phi_h)=l_h[f,P_hu](\phi_h^m),\quad\forall \phi_h^m\in\Phi_h^m.$$
Then
\begin{equation}\label{eqdgh1}
\|Q_h\psi-P_h\psi\|_{\mathrm{DG}}
\leq Ch^{-\frac12}\|P_hu-u\|_{L^2(\Gamma)}.
\end{equation}
Since $r=m-1$, Lemma~\ref{lem3} gives
\[
\|u-P_hu\|_{L^2(\Gamma)}\leq Ch^m,\quad
\|u-P_hu\|_U\leq Ch^{m-1}.
\]
It follows that
\begin{equation}\label{eq:Q-state-error}
\|Q_h\psi-\psi\|_{L^2(\Omega)}
\leq \|\psi-P_h\psi\|_{L^2(\Omega)}+Ch^{-\frac12}\|P_hu-u\|_{L^2(\Gamma)}
\leq Ch^{m-1}.
\end{equation}

Testing the continuous and discrete first order optimality conditions with $u_h^r-P_hu\in U_h^r$
gives
\[
\alpha\|u_h^r-P_hu\|_U^2=-l_h[0,u_h^r-P_hu](\eta_h^m-\eta).
\]
We can also show that $a_h(\psi_h^m-Q_h\psi,\phi_h)=l_h[0,u_h^r-P_hu](\phi_h)$. Taking
$\phi_h=\eta_h^m-P_h\eta$ and using \eqref{state1} yields
$$l_h[0,u_h^r-P_hu](\eta_h^m-P_h\eta)=(\psi_h^m-\psi,\psi_h^m-Q_h\psi).$$
Thus, the exact cancellation identity is
\begin{equation}\label{eq:control-cancellation}
\alpha\|u_h^r-P_hu\|_U^2+\|\psi_h^m-Q_h\psi\|_{L^2(\Omega)}^2
=-(Q_h\psi-\psi,\psi_h^m-Q_h\psi)-l_h[0,u_h^r-P_hu](P_h\eta-\eta).
\end{equation}
We note that this identity uses only discrete differences and the trace
$\Delta\eta|_\Gamma$; it remains valid for $\eta\in H^3(\Omega)$ and
an energy-regular state.

Combining Young's inequality, $\|u_h^r-P_hu\|_{L^2(\Gamma)}\leq\|u_h^r-P_hu\|_U$,
\eqref{eq:adjoint-boundary-projection} and
\eqref{eq:Q-state-error} we can obtain
\[
\|u_h^r-P_hu\|_U+\|\psi_h^m-Q_h\psi\|_{L^2(\Omega)}
\leq C\bigl(\|Q_h\psi-\psi\|_{L^2(\Omega)}
+h^{m-\frac32}\bigr)\leq Ch^{m-\frac32}.
\]
Adding the projection errors gives
$\|u-u_h^r\|_U+\|\psi-\psi_h^m\|_{L^2(\Omega)}
\leq Ch^{m-\frac32}$. Finally, \eqref{eqdg} and
\eqref{eq:projection-energy} give the two asserted DG estimates.
\end{proof}

\begin{thm}\label{thm2mod}
Under the hypotheses of Theorem~\ref{thm2}, with $m\geq3$,
\[
\interleave\psi-\psi_h^m\interleave_{\mathrm{DG}}\leq Ch^{m-2},
\quad
\interleave\eta-\eta_h^m\interleave_{\mathrm{DG}}\leq Ch^{m-\frac32}.
\]
\end{thm}
\begin{proof}
Now $\psi,\eta\in H^{m+1}(\Omega)\subseteq H^4(\Omega)$, so
strong DG norms exist. For $q=\psi,\eta$, with its corresponding
discrete approximation $q_h$, reconstruction and discrete norm
equivalence give
\[
\begin{aligned}
\interleave q-q_h\interleave_{\mathrm{DG}}
&\leq\interleave q-\mathcal R^mq\interleave_{\mathrm{DG}}
+C\|\mathcal R^mq-q_h\|_{\mathrm{DG}}\\
&\leq Ch^{m-1}\|q\|_{H^{m+1}(\Omega)}
+C\|q-q_h\|_{\mathrm{DG}}.
\end{aligned}
\]
Applying Theorem~\ref{thm2} gives the claimed result.
\end{proof}

\begin{rmk}\label{rmk:lowreg-control}
The smooth-solution hypotheses of Theorem~\ref{thm2} are
additional assumptions, not consequences of
Theorem~\ref{thm:regularity} for arbitrary degree.
The cancellation \eqref{eq:control-cancellation} also applies
with the actual low-regularity projection errors. On a general
convex polygon, however, the $h^{-\frac12}$ boundary estimate
\eqref{eqdgh1} combined with only
$u\in H^{\frac12}(\Gamma)$ does not yield a positive algebraic
convergence rate for the control variable. Under the angle restriction we instead use the
uniform $L^2$ stability proved below.
\end{rmk}

Let $\psi\in H^{s}(\Omega)$ for some $2\leq s\leq m+1$, we obtain in Theorem \ref{thm1} the convergence rate $O(h^{s-2})$ for the finite element approximation error to \eqref{eq8} under DG norm. This result is later used to prove error estimates for optimal control problems, yielding the convergence rate $O(h^{m-\frac{3}{2}})$ by assuming $H^{m+1}(\Omega)$ regularity for the solution. However, for biharmonic equations posed on polygonal domains, we can expect only at most $H^4(\Omega)$ regularity. This gives limited regularity for the optimal control, as stated in Theorem \ref{thm:regularity}. So, for practical optimal control problems posed on polygonal domains, we can expect only limited convergence rates. However, for manufactured solutions as constructed in Section 5, we can observe higher order convergence.

\subsection{Estimates under the available polygonal regularity}
In this subsection, we provide the error estimate for optimal control problems according to the regularity stated in Theorem \ref{thm:regularity}. We now impose the angle condition, but we have no additional smoothness
of the exact state.  

The following stability estimate for the discrete control-to-state mapping
allows us to use a smooth clamped adjoint for duality.

\begin{lem}\label{lem:state-L2-stability}
Let $m\geq3$ and suppose that all interior angles of $\Omega$
are less than $126.283$ degrees. Define
$S_h:L^2(\Gamma)\to\Phi_h^m$ by
\[
a_h(S_hv,\phi_h^m)=l_h[0,v](\phi_h^m),
\quad\forall \phi_h^m\in\Phi_h^m.
\]
Then
\begin{equation}\label{eq:state-L2-stability}
\|S_hv\|_{L^2(\Omega)}\leq C\|v\|_{L^2(\Gamma)}.
\end{equation}
\end{lem}
\begin{proof}
For $q\in L^2(\Omega)$, let $z\in H_0^2(\Omega)$ and
$z_h^m\in\Phi_h^m$ solve respectively the following continuous and discrete problems
\[
\Delta^2z=q,\quad a_h(\phi_h^m,z_h^m)=(\phi_h^m,q),
\quad\forall \phi_h^m\in\Phi_h^m.
\]
 Regularity result for biharmonic equations with clamped boundary conditions and Theorem~\ref{thm11} with $t=4$ give respectively
\[
\|z\|_{H^4(\Omega)}\leq C\|q\|_{L^2(\Omega)},\quad
\interleave z-z_h^m\interleave_{\mathrm{DG}}
\leq Ch^2\|q\|_{L^2(\Omega)}.
\]
The discrete equations imply
$(S_hv,q)=a_h(S_hv,z_h^m)=l_h[0,v](z_h^m)$.
Since $\pn z=0$, the trace theorem gives
$|l_h[0,v](z)|=|(v,\Delta z)_\Gamma|
\leq C\|v\|_{L^2(\Gamma)}\|q\|_{L^2(\Omega)}$.
Moreover, \eqref{eq:boundary-functional-strong} yields
\[
|l_h[0,v](z_h^m-z)|
\leq Ch^{-\frac12}\|v\|_{L^2(\Gamma)}
\interleave z_h^m-z\interleave_{\mathrm{DG}}
\leq Ch^{\frac32}\|v\|_{L^2(\Gamma)}\|q\|_{L^2(\Omega)}.
\]
Therefore, 
\begin{eqnarray}
\|S_hv\|_{L^2(\Omega)}&=&\sup\limits_{q\in L^2(\Omega),q\neq 0}\frac{(S_hv,q)}{\|q\|_{L^2(\Omega)}}=\sup\limits_{q\in L^2(\Omega),q\neq 0}\frac{l_h[0,v](z_h^m)}{\|q\|_{L^2(\Omega)}}\nonumber\\
&\leq&\sup\limits_{q\in L^2(\Omega),q\neq 0}\frac{|l_h[0,v](z)|}{\|q\|_{L^2(\Omega)}}+\sup\limits_{q\in L^2(\Omega),q\neq 0}\frac{|l_h[0,v](z-z_h^m)|}{\|q\|_{L^2(\Omega)}}\nonumber\\
&\leq &C(h^{\frac32}\|v\|_{L^2(\Gamma)}+\|v\|_{L^2(\Gamma)}).\nonumber
\end{eqnarray}
This proves the result.
\end{proof}

\begin{thm}\label{thm3mod}
Let $f,\psi_d\in L^2(\Omega)$ and let the continuous and discrete
optimal solutions satisfy \eqref{eq6} and \eqref{eq13}.
Take $r=m-1$, $m\geq3$, and suppose that all interior angles
of $\Omega$ are less than $126.283$ degrees.
If $U=L^2(\Gamma)$, then
\begin{equation}\label{eq:lowreg-L2-rates}
\|u-u_h^r\|_{L^2(\Gamma)}
+\|\psi-\psi_h^m\|_{L^2(\Omega)}
+\|\eta-\eta_h^m\|_{L^2(\Omega)}
\leq Ch^{\frac32}.
\end{equation}
If $U=H^1(\Gamma)$, then
\begin{equation}\label{eq:lowreg-H1-rates}
\begin{aligned}
\|u-u_h^r\|_{H^1(\Gamma)}&\leq C_\delta h^\frac32,\\
\|u-u_h^r\|_{L^2(\Gamma)}
+\|\psi-\psi_h^m\|_{L^2(\Omega)}
+\|\eta-\eta_h^m\|_{L^2(\Omega)}&\leq C_\delta h^{\frac32}.
\end{aligned}
\end{equation}
In the second case any fixed $0<\delta<1/2$ may be used in
the regularity estimate for $\psi$; the constant is not asserted
to be uniform as $\delta\to0$. The constants depend on the
data and the corresponding norms in Theorem~\ref{thm:regularity}.
\end{thm}
\begin{proof}
Let $s=3$ if $U=L^2(\Gamma)$, and fix $s=3-\delta$,
$0<\delta<1/2$, if $U=H^1(\Gamma)$. By
Theorem~\ref{thm:regularity},
$\psi\in H^s(\Omega)$ and $\eta\in H^4(\Omega)\cap H_0^2(\Omega)$.
We use the projections in \eqref{eq:elliptic-projections}.
The $L^2$ norm error estimate of Theorem~\ref{thm1}, with $\epsilon=2$,
implies
\begin{equation}\label{eq:lowreg-state-projection}
\|\psi-P_h\psi\|_{L^2(\Omega)}
\leq C\left(h^s\|\psi\|_{H^s(\Omega)}
+h^2\operatorname{Osc}(f)\right)
\leq C(h^s+h^4).
\end{equation}
Since $\eta\in H^4(\Omega)$, Theorem~\ref{thm11} with $t=4$ and
\eqref{eq:boundary-functional-strong} give
\begin{equation}\label{eq:lowreg-adjoint-projection}
\begin{aligned}
\interleave\eta-P_h\eta\interleave_{\mathrm{DG}}
\leq Ch^2\|\eta\|_{H^4(\Omega)},\quad 
|l_h[0,v](P_h\eta-\eta)|
\leq Ch^{\frac32}\|v\|_{L^2(\Gamma)},
\quad\forall v\in L^2(\Gamma).
\end{aligned}
\end{equation}
We remark that only the adjoint error, not the state error, is measured in
the strong DG norm.

Let $P_hu$ be the control projection from Lemma~\ref{lem3}. In the $L^2$ case the trace theorem gives
$u\in H_{\rm pw}^{\frac32}(\Gamma)$; in the $H^1$ case
$u\in H_{\rm pw}^2(\Gamma)\cap U_{ad}$. Remark
\ref{cor:control-projection}, applicable since $r\geq2$, yields
\begin{equation}\label{eq:lowreg-control-projection}
\begin{array}{lll}
U=L^2(\Gamma):&\|u-P_hu\|_{L^2(\Gamma)}\leq Ch^{\frac32},\\
U=H^1(\Gamma):&\|u-P_hu\|_{L^2(\Gamma)}\leq Ch^2,\quad \|u-P_hu\|_{H^1(\Gamma)}\leq Ch^{s-1}.
\end{array}
\end{equation}
Define $Q_h\psi$ as in the proof of Theorem~\ref{thm2}.
The discrete equations give
$Q_h\psi-P_h\psi=S_h(P_hu-u)$, so
Lemma~\ref{lem:state-L2-stability} implies
\begin{equation}\label{eq:lowreg-Q-error}
\|Q_h\psi-\psi\|_{L^2(\Omega)}\leq \|\psi-P_h\psi\|_{L^2(\Omega)}+C\|u-P_hu\|_{L^2(\Gamma)}.
\end{equation}
This is an $L^2$ norm error estimate without factor $h^{-\frac12}$, compared to the proof of Theorem \ref{thm2}.

Note that $\psi_h^m-Q_h\psi=S_h(u_h^r-P_hu)$.
The algebraic identity \eqref{eq:control-cancellation} is
available at the present regularity. Explicitly,
$$\alpha\|u_h^r-P_hu\|_U^2=-l_h[0,u_h^r-P_hu](\eta_h^m-\eta),\quad l_h[0,u_h^r-P_hu](\eta_h^m-P_h\eta)=(\psi_h^m-\psi,\psi_h^m-Q_h\psi),$$ hence
\[
\alpha\|u_h^r-P_hu\|_U^2+\|\psi_h^m-Q_h\psi\|_{L^2(\Omega)}^2
=-(Q_h\psi-\psi,\psi_h^m-Q_h\psi)-l_h[0,u_h^r-P_hu](P_h\eta-\eta).
\]
By \eqref{eq:lowreg-adjoint-projection}, Young's inequality,
and \eqref{eq:lowreg-Q-error}, we can conclude
\begin{equation}\label{eq:lowreg-superclose}
\|u_h^r-P_hu\|_U+\|\psi_h^m-Q_h\psi\|_{L^2(\Omega)}
\leq C\bigl(\|\psi-P_h\psi\|_{L^2(\Omega)}+\|u-P_hu\|_{L^2(\Gamma)}+h^{\frac32}\bigr).
\end{equation}
Therefore
\begin{equation}\label{eq:lowreg-control-state}
\begin{aligned}
\|u-u_h^r\|_U-\|u-P_hu\|_U&\leq
C\bigl(\|\psi-P_h\psi\|_{L^2(\Omega)}+\|u-P_hu\|_{L^2(\Gamma)}+h^{\frac32}\bigr),\\
\|u-u_h^r\|_{L^2(\Gamma)}
+\|\psi-\psi_h^m\|_{L^2(\Omega)}
&\leq C\bigl(\|\psi-P_h\psi\|_{L^2(\Omega)}+\|u-P_hu\|_{L^2(\Gamma)}+h^{\frac32}\bigr).
\end{aligned}
\end{equation}
For the estimate of the second line we used $\|u_h^r-P_hu\|_{L^2(\Gamma)}\leq\|u_h^r-P_hu\|_U$.
Since $s>5/2$, equations
\eqref{eq:lowreg-state-projection} and
\eqref{eq:lowreg-control-projection} prove the control and
state bounds in the statement.

Finally subtract the adjoint equations in \eqref{eq13} and
\eqref{eq:elliptic-projections}. Coercivity of $a_h$ and Lemma~\ref{lem1}
give
\[
\|\eta_h^m-P_h\eta\|_{\mathrm{DG}}
\leq C\|\psi_h^m-\psi\|_{L^2(\Omega)}.
\]
Using Lemma~\ref{lem1} once more and
\eqref{eq:lowreg-adjoint-projection}, we obtain
\[
\|\eta-\eta_h^m\|_{L^2(\Omega)}
\leq Ch^2\|\eta\|_{H^4(\Omega)}
+C\|\psi_h^m-\psi\|_{L^2(\Omega)}\leq Ch^{\frac32}.
\]
This proves both cases without a strong DG norm or Galerkin
orthogonality for the state.
\end{proof}

\begin{rmk}\label{rmk:lowreg-energy}
The same argument also gives energy estimates at the actual
regularity. Theorem~\ref{thm1} and \eqref{eqdg} imply
\[
\|\psi-\psi_h^m\|_{\mathrm{DG}}
\leq C\bigl(h^{s-2}+h^2+
h^{-\frac12}\|u-u_h^r\|_{L^2(\Gamma)}\bigr).
\]
Thus the convergence rate for the energy norm error of state is $O(h)$ for $U=L^2(\Gamma)$ and
$O(h^{1-\delta})$ for $U=H^1(\Gamma)$. For the adjoint variable, add and
subtract $P_h\eta$, use discrete norm equivalence and
\eqref{eq:lowreg-adjoint-projection}, we can  obtain
$\interleave\eta-\eta_h^m\interleave_{\mathrm{DG}}\leq Ch^{\frac32}$.
The state remains measured only in the weak DG norm.
\end{rmk}

\section{Numerical Results}\label{sec4}
Now we are going to demonstrate the theoretical results 
in Theorem \ref{thm2mod},
Theorem \ref{thm3mod} and Remark \ref{rmk1}, 
through three examples, with 
$U=L^2(\Gamma)$ and $H^1(\Gamma)$. Recall that  
\[e_\psi:=\psi^*-\psi_h^{m*},\quad e_u:=u^*-u_h^{r*},\quad e_\eta:=\eta^*-\eta_h^{m*}.\]
The columns of order in the following tables are calculated by $\log_2\frac
{\mathrm{error}_{\rm old}}{\mathrm{error}_{\rm new}}$. We all take $r=m-1$ in the following
examples.

\begin{exe}\label{exe1}
Take $\Omega=(0,1)^2$, $U=L^2(\Gamma)$ and $\alpha=1$. Consider the following optimization problem
\begin{equation*}
	\begin{aligned}
		&\min_{u\in U}\frac{1}{2}\|\psi-\psi_d\|_{L^2(\Omega)}^2+\frac12
		\|u\|_{L^2(\Gamma)}^2,&\\
		\mathrm{s.t.}\ 
		&\Delta^2 \psi= f\quad\mathrm{in}\ 
		\Omega,\\
		&\psi=0,\ \pn\psi 
		=u\quad \mathrm{on}\ \Gamma,
	\end{aligned}
\end{equation*}
where 
\begin{equation*}
	\begin{aligned}
		&\psi=-2{\pi^2}[x_1(1-x_1)\mathrm{sin}^2({\pi x_2})
		+x_2(1-x_2)\mathrm{sin}^2({\pi x_1})],&\\
		&\eta={}\mathrm{sin}^2({\pi x_1})\mathrm{sin}^2({\pi x_2}),
		\quad u=\pn \psi|_\Gamma,\\
		&f=\Delta^2\psi,\quad \psi_d=\psi-\Delta^2\eta.
	\end{aligned}
\end{equation*}
The exact solutions are given by  
\[\psi^*=\psi,\quad \eta^*=\eta,\quad u^*=u.\]
\end{exe}
Here we use the classical 
$r$-th order Lagrange element on each $\Gamma_i$ ($i=1,2,3,4$), denoted as $U_{i,h}^r$.
Define $U_{h} := \ds\prod_{i=1}^l U_{i,h}^r.$ Take $m=2,3,4$, the 
errors are shown below. Here we take $\mu=2.5 m^2$. 
\begin{table}[H]
	\centering
	\renewcommand\arraystretch{1.0}
	\begin{tabular}{l|l|l|l|l|l|l|l|l}
		\hline\hline
		$h$&$\|e_{\psi}\|_{\mathrm{DG}}$&order&$\|e_{\eta}\|_{\mathrm{DG}}$&order&$|e_{\psi}|_{H^1(\Omega)}$&order&$|e_\eta|_{H^1(\Omega)}$&order\\
		\hline
		1/8&3.66e1&&3.70&&2.45&&2.22e-1&\\
		\hline
		1/16&1.48e1&1.31&1.65&1.16&4.06e-1&2.59&4.17e-2&2.41\\
		\hline
		1/32&6.29&1.23&8.16e-1&1.02&8.08e-2&2.33&1.00e-2&2.06\\
		\hline
		1/64&2.85&1.14&4.05e-1&1.01&1.84e-2&2.13&2.38e-3&2.07\\
		\hline
		1/128&1.40&1.03&2.02e-1&1.00&4.51e-3&2.03&5.72e-4&2.06\\
		\hline\hline
	\end{tabular}\vskip-0.3cm
	\caption{DG and $H^1$ norms of $e_\psi$ and $e_\eta$ in Example 1 for $m=2$.}\label{tab1}
\end{table}

\begin{table}[H]
	\centering
	\renewcommand\arraystretch{1.0}
	\begin{tabular}{l|l|l|l|l|l|l}
		\hline\hline
		$h$&$\|e_{\psi}\|_{L^2(\Omega)}$&order&$\|e_{\eta}\|_{L^2(\Omega)}$&order&$\|e_{u}\|_{L^2(\Gamma)}$&order\\
		\hline
		1/8&5.52e-1&&2.26e-2&&1.29e1&\\
		\hline
		1/16&9.00e-2&2.62&1.96e-3&3.53&3.15&2.03\\
		\hline
		1/32&1.86e-2&2.27&3.64e-4&2.43&7.47e-1&2.08\\
		\hline
		1/64&4.14e-3&2.17&8.57e-5&2.09&1.93e-1&1.95\\
		\hline
		1/128&9.58e-4&2.11&2.37e-5&1.85&6.27e-2&1.62\\
		\hline\hline
	\end{tabular}\vskip-0.3cm
	\caption{$L^2$ norms of $e_\psi,e_\eta$ and $e_u$ in Example 1 for $m=2$.}
\end{table}


\begin{table}[H]
	\centering
	\renewcommand\arraystretch{1.0}
	\begin{tabular}{l|l|l|l|l|l|l|l|l}
		\hline\hline
		$h$&$\|e_{\psi}\|_{\mathrm{DG}}$&order&$\|e_{\eta}\|_{\mathrm{DG}}$&order&$|e_{\psi}|_{H^1(\Omega)}$&order&$|e_\eta|_{H^1(\Omega)}$&order\\
		\hline
		1/8&1.97e1&&1.80&&1.13&&8.36e-2&\\
		\hline
		1/16&9.36&1.07&4.20e-1&2.10&2.90e-1&1.96&8.73e-3&3.26\\
		\hline
		1/32&2.83&1.73&9.71e-2&2.11&4.27e-2&2.76&9.03e-4&3.27\\
		\hline
		1/64&7.41e-1&1.93&2.33e-2&2.06&5.55e-3&2.94&9.85e-5&3.20\\
		\hline
		1/128&1.87e-1&1.99&5.71e-3&2.03&7.02e-4&2.98&1.12e-5&3.14\\
		\hline\hline
	\end{tabular}\vskip-0.3cm
	\caption{DG and $H^1$ norms of $e_\psi$ and $e_\eta$ in Example 1 for $m=3$.}
\end{table}

\begin{table}[H]
	\centering
	\renewcommand\arraystretch{1.0}
	\begin{tabular}{l|l|l|l|l|l|l}
		\hline\hline
		$h$&$\|e_{\psi}\|_{L^2(\Omega)}$&order&$\|e_{\eta}\|_{L^2(\Omega)}$&order&$\|e_{u}\|_{L^2(\Gamma)}$&order\\
		\hline
		1/8&8.44e-2&&1.25e-2&&6.70&\\
		\hline
		1/16&1.78e-2&2.24&7.84e-4&4.00&2.44&1.46\\
		\hline
		1/32&1.27e-3&3.81&4.89e-5&4.00&5.31e-1&2.20\\
		\hline
		1/64&7.34e-5&4.11&3.13e-6&3.97&1.13e-1&2.23\\
		\hline
		1/128&4.32e-6&4.09&2.03e-7&3.95&2.48e-2&2.19\\
		\hline\hline
	\end{tabular}\vskip-0.3cm
	\caption{$L^2$ norms of $e_\psi,e_\eta$ and $e_u$ in Example 1 for $m=3$.}
\end{table}


\begin{table}[H]
	\centering
	\renewcommand\arraystretch{1.0}
	\begin{tabular}{l|l|l|l|l|l|l|l|l}
		\hline\hline
		$h$&$\|e_{\psi}\|_{\mathrm{DG}}$&order&$\|e_{\eta}\|_{\mathrm{DG}}$&order&$|e_{\psi}|_{H^1(\Omega)}$&order&$|e_\eta|_{H^1(\Omega)}$&order\\
		\hline
		1/8&5.64e1&&7.16e-1&&2.85&&3.83e-2&\\
		\hline
		1/16&1.26e1&2.16&8.22e-2&3.12&2.66e-1&3.42&1.75e-3&4.45\\
		\hline
		1/32&2.28&2.47&1.00e-2&3.04&2.33e-2&3.51&1.01e-4&4.11\\
		\hline
		1/64&4.14e-1&2.46&1.12e-3&3.16&2.11e-3&3.47&5.73e-6&4.14\\
		\hline
		1/128&5.54e-2&2.90&1.26e-4&3.15&1.40e-4&3.91&3.28e-7&4.13\\
		\hline\hline
	\end{tabular}\vskip-0.3cm
	\caption{DG and $H^1$ norms of $e_\psi$ and $e_\eta$ in Example 1 for $m=4$.}
\end{table}

\begin{table}[H]
	\centering
	\renewcommand\arraystretch{1.0}
	\begin{tabular}{l|l|l|l|l|l|l}
		\hline\hline
		$h$&$\|e_{\psi}\|_{L^2(\Omega)}$&order&$\|e_{\eta}\|_{L^2(\Omega)}$&order&$\|e_{u}\|_{L^2(\Gamma)}$&order\\
		\hline
		1/8&3.52e-1&&4.13e-3&&1.04e1&\\
		\hline
		1/16&1.17e-2&4.91&8.37e-5&5.62&1.32&2.98\\
		\hline
		1/32&3.62e-4&5.01&1.51e-6&5.79&1.90e-1&2.80\\
		\hline
		1/64&1.46e-5&4.63&3.20e-8&5.56&2.73e-2&2.80\\
		\hline
		1/128&4.75e-7&4.94&7.6e-10&5.40&3.67e-3&2.90\\
		\hline\hline
	\end{tabular}\vskip-0.3cm
	\caption{$L^2$ norms of $e_\psi,e_\eta$ and $e_u$ in Example 1 for $m=4$.}\label{tab6}
\end{table}

From Tables \ref{tab1}-\ref{tab6} we can see that the convergence orders of 
the energy norm errors are close to $m-1$, while those of the $H^1$ 
semi-norms and $L^2$ norms are  more than $m$ and even close to $m+1$. These results are 
generally better than the theoretical estimates provided in the 
above theorems. For the convergence of control, we can observe
an approximately $(m-1)$-th order convergence, which is also sharper 
than the value $m-\ds\frac32$ given in Theorem \ref{thm2}.

\begin{exe}\label{exe2}
	Take $\Omega=(0,1)^2$, $U=H^1(\Gamma)$ and $\alpha=1$. Consider the following optimization problem
	\begin{equation*}
		\begin{aligned}
			&\min_{u\in U_{ad}}\frac{1}{2}\|\psi-\psi_d\|_{L^2(\Omega)}^2+\frac12
			\|u\|_{H^1(\Gamma)}^2,&\\
			\mathrm{s.t.}\ 
			&\Delta^2 \psi= f\quad\mathrm{in}\ 
			\Omega,\\
			&\psi=0,\ \pn\psi 
			=u\quad \mathrm{on}\ \Gamma,
		\end{aligned}
	\end{equation*}
	where 
	\begin{equation*}
		\begin{aligned}
			&\psi=-\frac{1}{\pi}[\mathrm{sin}^4(\pi x_1)
			\mathrm{sin}(\pi x_2)+
			\mathrm{sin}^4(\pi x_2)
			\mathrm{sin}(\pi x_1)],&\\
			&\eta=\frac12\mathrm{sin}^2(\pi x_1)\mathrm{sin}^2(\pi x_2)
			\Big[\Big(\frac{1}{\pi^2}+16\Big)
			(\mathrm{sin}^2(\pi x_1)+\mathrm{sin}^2(\pi x_2))-12\Big],
			\\
			&f=\Delta^2\psi,\quad \psi_d=\psi-\Delta^2\eta,\quad u=\pn \psi|_\Gamma.
		\end{aligned}
	\end{equation*}
The exact solutions are given by 
\[\psi^*=\psi,\quad \eta^*=\eta,\quad u^*=u.\]
\end{exe}
%
Here we use a subspace of   
$r$-th order Lagrange element on each $\Gamma_i$ with the value on its endpoints forced to be 0,
denoted as $U_{i,h}^r$.
Define $U_{h} := \ds\prod_{i=1}^l U_{i,h}^r.$
Take $m=2,3,4$, the 
errors are shown below. Here we take $\mu=2.5m^2.$

\begin{table}[H]
	\centering
	\renewcommand\arraystretch{1.0}
	\begin{tabular}{l|l|l|l|l|l|l|l|l}
		\hline\hline
		$h$&$\|e_{\psi}\|_{\mathrm{DG}}$&order&$\|e_{\eta}\|_{\mathrm{DG}}$&order&$|e_{\psi}|_{H^1(\Omega)}$&order&$|e_\eta|_{H^1(\Omega)}$&order\\
		\hline
		1/8&4.22&&5.13e1&&3.02e-1&&1.69&\\
		\hline
		1/16&1.78&1.25&2.58e1&0.99&8.26e-2&1.87&4.46e-1&1.92\\
		\hline
		1/32&7.10e-1&1.33&1.30e1&0.99&2.12e-2&1.96&1.17e-1&1.93\\
		\hline
		1/64&3.07e-1&1.21&6.49&1.00&5.02e-3&2.08&2.97e-2&1.98\\
		\hline
		1/128&1.43e-1&1.10&3.25&1.00&1.19e-3&2.08&7.38e-3&2.01\\
		\hline\hline
	\end{tabular}\vskip-0.3cm
	\caption{DG and $H^1$ norms of $e_\psi$ and $e_\eta$ in Example 2 for $m=2$.}\label{tb7}
\end{table}

\begin{table}[H]
	\centering
	\renewcommand\arraystretch{1.0}
	\begin{tabular}{l|l|l|l|l|l|l}
		\hline\hline
		$h$&$\|e_{\psi}\|_{L^2(\Omega)}$&order&$\|e_{\eta}\|_{L^2(\Omega)}$&order&$\|e_{u}\|_{H^1(\Gamma)}$&order\\
		\hline
		1/8&6.67e-2&&1.14e-1&&4.26&\\
		\hline
		1/16&1.72e-2&1.96&2.66e-2&2.10&1.82&1.23\\
		\hline
		1/32&4.42e-3&1.96&6.46e-3&2.04&5.47e-1&1.73\\
		\hline
		1/64&1.02e-3&2.12&1.61e-3&2.00&2.02e-1&1.44\\
		\hline
		1/128&2.42e-4&2.08&4.11e-4&1.97&9.18e-2&1.14\\
		\hline\hline
	\end{tabular}\vskip-0.3cm
	\caption{$L^2$ norms of $e_\psi,e_\eta$ and $H^1$ norm of $e_u$ in Example 2 for $m=2$.}
\end{table}


\begin{table}[H]
	\centering
	\renewcommand\arraystretch{1.0}
	\begin{tabular}{l|l|l|l|l|l|l|l|l}
		\hline\hline
		$h$&$\|e_{\psi}\|_{\mathrm{DG}}$&order&$\|e_{\eta}\|_{\mathrm{DG}}$&order&$|e_{\psi}|_{H^1(\Omega)}$&order&$|e_\eta|_{H^1(\Omega)}$&order\\
		\hline
		1/8&2.43&&4.27e1&&1.29e-1&&2.17&\\
		\hline
		1/16&4.03e-1&2.59&9.92&2.11&9.71e-3&3.73&2.29e-1&3.24\\
		\hline
		1/32&9.38e-2&2.10&2.35&2.08&9.73e-4&3.32&2.37e-2&3.27\\
		\hline
		1/64&2.24e-2&2.07&5.69e-1&2.05&1.02e-4&3.25&2.53e-3&3.23\\
		\hline
		1/128&5.49e-3&2.03&1.41e-1&2.01&1.14e-5&3.16&2.92e-4&3.12\\
		\hline\hline
	\end{tabular}\vskip-0.3cm
	\caption{DG and $H^1$ norms of $e_\psi$ and $e_\eta$ in Example 2 for $m=3$.}
\end{table}

\begin{table}[H]
	\centering
	\renewcommand\arraystretch{1.0}
	\begin{tabular}{l|l|l|l|l|l|l}
		\hline\hline
		$h$&$\|e_{\psi}\|_{L^2(\Omega)}$&order&$\|e_{\eta}\|_{L^2(\Omega)}$&order&$\|e_{u}\|_{H^1(\Gamma)}$&order\\
		\hline
		1/8&1.15e-2&&1.86e-1&&5.05&\\
		\hline
		1/16&6.69e-4&4.10&1.43e-2&3.70&3.03e-1&4.06\\
		\hline
		1/32&3.79e-5&4.14&1.07e-3&3.74&3.74e-2&3.02\\
		\hline
		1/64&2.98e-6&3.67&7.60e-5&3.82&5.08e-3&2.88\\
		\hline
		1/128&2.28e-7&3.71&5.13e-6&3.89&9.61e-4&2.40\\
		\hline\hline
	\end{tabular}\vskip-0.3cm
	\caption{$L^2$ norms of $e_\psi,e_\eta$ and $H^1$ norm of $e_u$ in Example 2 for $m=3$.}
\end{table}


\begin{table}[H]
	\centering
	\renewcommand\arraystretch{1.0}
	\begin{tabular}{l|l|l|l|l|l|l|l|l}
		\hline\hline
		$h$&$\|e_{\psi}\|_{\mathrm{DG}}$&order&$\|e_{\eta}\|_{\mathrm{DG}}$&order&$|e_{\psi}|_{H^1(\Omega)}$&order&$|e_\eta|_{H^1(\Omega)}$&order\\
		\hline
		1/8&2.93&&3.47e1&&1.56e-1&&1.40&\\
		\hline
		1/16&2.02e-1&3.86&4.01&3.11&4.82e-3&5.02&8.27e-2&4.08\\
		\hline
		1/32&1.46e-2&3.79&3.98e-1&3.33&1.44e-4&5.06&4.17e-3&4.31\\
		\hline
		1/64&1.55e-3&3.24&3.87e-2&3.36&7.44e-6&4.27&2.06e-4&4.34\\
		\hline
		1/128&1.81e-4&3.10&4.09e-3&3.24&4.63e-7&4.01&1.12e-5&4.20\\
		\hline\hline
	\end{tabular}\vskip-0.3cm
	\caption{DG and $H^1$ norms of $e_\psi$ and $e_\eta$ in Example 2 for $m=4$.}
\end{table}

\begin{table}[H]
	\centering
	\renewcommand\arraystretch{1.0}
	\begin{tabular}{l|l|l|l|l|l|l}
		\hline\hline
		$h$&$\|e_{\psi}\|_{L^2(\Omega)}$&order&$\|e_{\eta}\|_{L^2(\Omega)}$&order&$\|e_{u}\|_{H^1(\Gamma)}$&order\\
		\hline
		1/8&1.68e-2&&8.29e-2&&8.31&\\
		\hline
		1/16&3.03e-4&5.79&2.21e-3&5.23&4.79e-1&4.12\\
		\hline
		1/32&7.31e-6&5.37&5.21e-5&5.41&2.36e-2&4.34\\
		\hline
		1/64&1.24e-7&5.88&1.10e-6&5.57&1.76e-3&3.75\\
		\hline
		1/128&2.31e-9&5.75&2.34e-8&5.55&1.30e-4&3.76\\
		\hline\hline
	\end{tabular}\vskip-0.3cm
	\caption{$L^2$ norms of $e_\psi,e_\eta$ and $H^1$ norm of $e_u$ in Example 2 for $m=4$.}\label{tb12}
\end{table}

Again, from Table \ref{tb7}-\ref{tb12} we can see that, the convergence orders of 
the energy norms are no less than $m-1$, while those of the $H^1$ 
semi-norms and $L^2$ norms are near to $m$. These results are 
also better than the theoretical estimates proved in the 
above theorems. For the convergence of control, we can observe
an approximately first order convergence for $m=2$ and generally 
$m$-th order 
convergence for $m=3,4$. This is also better 
than the results in Theorem \ref{thm2}.

Next, we consider the case without exact solutions, where $\Omega$ is set as a rhombus. Just like Example \ref{exe1}, we take $\mu=2.5 m^2$.

\begin{exe}\label{exe3}
	Take $\Omega\subseteq\mathbb{R}^2$, $U=L^2(\Gamma)$ and $\alpha=1$. Consider the following optimal control problem
	\begin{equation*}
		\begin{aligned}
			&\min_{u\in U_{ad}}\frac{1}{2}\|\psi-\psi_d\|_{L^2(\Omega)}^2+\frac12
			\|u\|_{L^2(\Gamma)}^2,&\\
			\mathrm{s.t.}\ 
			&\Delta^2 \psi= f\quad\mathrm{in}\ 
			\Omega,\\
			&\psi=0,\ \pn\psi 
			=u\quad \mathrm{on}\ \Gamma,
		\end{aligned}
	\end{equation*}
	where 
	$\Omega$ is the domain circled by the following four lines:
	\[x_2-\frac25x_1=0,\ x_2+\frac25x_1=0,\ x_2-\frac25x_1+\frac45=0,\ x_2+\frac25x_1-\frac45=0,\] and
	\begin{equation*}
		\begin{aligned}
			&f=1+0.3x_1-0.2x_2+0.1x_1x_2,\\&\psi_d= \sin\Big(\pi\Big(\frac{x_1}{2}+\frac{5x_2}{4}\Big)\Big)+
			0.37\cos\Big(2\pi\Big(\frac{x_1}{2}-\frac{5x_2}{4}\Big)\Big)+
			0.21\Big(\frac{x_1}{2}+\frac{5x_2}{4}\Big)\Big(\frac{x_1}{2}-\frac{5x_2}{4}\Big).
		\end{aligned}
	\end{equation*}
\end{exe}

Since there are no exact solutions, we define \[
\widetilde{e}_\psi(h):=\psi_h^{m*}-\psi_{h/2}^{m*},\quad
\widetilde{e}_\eta(h):=\eta_h^{m*}-\eta_{h/2}^{m*},\quad
\widetilde{e}_u(h):=u_h^{r*}-u_{h/2}^{r*}\]
and compare errors between two adjacent refinement levels.

\begin{table}[H]
	\centering
	\renewcommand\arraystretch{1.0}
	\begin{tabular}{l|l|l|l|l|l|l|l|l}
		\hline\hline
		$h$&$\|\widetilde{e}_{\psi}\|_{\mathrm{DG}}$&order&$\|\widetilde{e}_{\eta}\|_{\mathrm{DG}}$&order&$|\widetilde{e}_{\psi}|_{H^1(\Omega)}$&order&$|\widetilde{e}_\eta|_{H^1(\Omega)}$&order\\
		\hline
		1/8&2.85e-2&&3.99e-3&&7.49e-4&&7.39e-5&\\
		\hline
		1/16&2.25e-2&0.34&1.83e-3&1.12&4.19e-4&0.84&1.69e-5&2.13\\
		\hline
		1/32&1.71e-2&0.40&8.76e-4&1.06&1.92e-4&1.13&4.50e-6&1.91\\
		\hline
		1/64&1.15e-2&0.57&4.25e-4&1.04&7.24e-5&1.41&1.21e-6&1.89\\
		\hline
		1/128&7.47e-3&0.62&2.09e-4&1.02&2.38e-5&1.61&3.08e-7&1.97\\
		\hline\hline
	\end{tabular}\vskip-0.3cm
	\caption{DG and $H^1$ norms of $\widetilde{e}_\psi$ and $\widetilde{e}_\eta$ in Example 3 for $m=2$.}\label{tb13}
\end{table}

\begin{table}[H]
	\centering
	\renewcommand\arraystretch{1.0}
	\begin{tabular}{l|l|l|l|l|l|l}
		\hline\hline
		$h$&$\|\widetilde{e}_{\psi}\|_{L^2(\Omega)}$&order&$\|\widetilde{e}_{\eta}\|_{L^2(\Omega)}$&order&$\|\widetilde{e}_{u}\|_{L^2(\Gamma)}$&order\\
		\hline
		1/8&6.08e-5&&1.66e-5&&2.82e-3&\\
		\hline
		1/16&4.87e-5&0.32&2.80e-6&2.57&1.56e-3&0.85\\
		\hline
		1/32&2.32e-5&1.07&4.77e-7&2.55&1.04e-3&0.58\\
		\hline
		1/64&7.95e-6&1.55&8.98e-8&2.41&5.96e-4&0.80\\
		\hline
		1/128&2.29e-6&1.80&2.03e-8&2.15&3.08e-4&0.95\\
		\hline\hline
	\end{tabular}\vskip-0.3cm
	\caption{$L^2$ norms of $\widetilde{e}_\psi,\widetilde{e}_\eta$ and $\widetilde{e}_u$ in Example 3 for $m=2$.}
\end{table}


\begin{table}[H]
	\centering
	\renewcommand\arraystretch{1.0}
	\begin{tabular}{l|l|l|l|l|l|l|l|l}
		\hline\hline
		$h$&$\|\widetilde{e}_{\psi}\|_{\mathrm{DG}}$&order&$\|\widetilde{e}_{\eta}\|_{\mathrm{DG}}$&order&$|\widetilde{e}_{\psi}|_{H^1(\Omega)}$&order&$|\widetilde{e}_\eta|_{H^1(\Omega)}$&order\\
		\hline
		1/8&7.59e-2&&1.94e-3&&3.42e-3&&8.21e-5&\\
		\hline
		1/16&6.64e-2&0.19&6.64e-4&1.55&1.86e-3&0.88&1.61e-5&2.35\\
		\hline
		1/32&4.12e-2&0.69&1.97e-4&1.75&5.87e-4&1.66&2.33e-6&2.79\\
		\hline
		1/64&2.31e-2&0.83&5.55e-5&1.83&1.68e-4&1.81&3.22e-7&2.86\\
		\hline
		1/128&1.21e-2&0.93&1.50e-5&1.89&4.49e-5&1.90&4.28e-8&2.91\\
		\hline\hline
	\end{tabular}\vskip-0.3cm
	\caption{DG and $H^1$ norms of $\widetilde{e}_\psi$ and $\widetilde{e}_\eta$ in Example 3 for $m=3$.}
\end{table}

\begin{table}[H]
	\centering
	\renewcommand\arraystretch{1.0}
	\begin{tabular}{l|l|l|l|l|l|l}
		\hline\hline
		$h$&$\|\widetilde{e}_{\psi}\|_{L^2(\Omega)}$&order&$\|\widetilde{e}_{\eta}\|_{L^2(\Omega)}$&order&$\|\widetilde{e}_{u}\|_{L^2(\Gamma)}$&order\\
		\hline
		1/8&3.51e-4&&8.18e-6&&1.13e-2&\\
		\hline
		1/16&1.78e-4&0.98&1.57e-6&2.38&7.63e-3&0.57\\
		\hline
		1/32&3.78e-5&2.24&1.78e-7&3.14&3.46e-3&1.14\\
		\hline
		1/64&6.67e-6&2.50&1.70e-8&3.39&1.43e-3&1.27\\
		\hline
		1/128&1.04e-6&2.68&1.48e-9&3.52&5.51e-4&1.38\\
		\hline\hline
	\end{tabular}\vskip-0.3cm
	\caption{$L^2$ norms of $\widetilde{e}_\psi,\widetilde{e}_\eta$ and $\widetilde{e}_u$ in Example 3 for $m=3$.}
\end{table}


\begin{table}[H]
	\centering
	\renewcommand\arraystretch{1.0}
	\begin{tabular}{l|l|l|l|l|l|l|l|l}
		\hline\hline
		$h$&$\|\widetilde{e}_{\psi}\|_{\mathrm{DG}}$&order&$\|\widetilde{e}_{\eta}\|_{\mathrm{DG}}$&order&$|\widetilde{e}_{\psi}|_{H^1(\Omega)}$&order&$|\widetilde{e}_\eta|_{H^1(\Omega)}$&order\\
		\hline
		1/8&4.29e-2&&7.29e-4&&1.74e-3&&2.65e-5&\\
		\hline
		1/16&2.62e-2&0.71&2.22e-4&1.72&6.22e-4&1.48&3.86e-6&2.78\\
		\hline
		1/32&1.22e-2&1.10&6.24e-5&1.83&1.35e-4&2.20&5.42e-7&2.83\\
		\hline
		1/64&5.60e-3&1.12&1.63e-5&1.94&2.94e-5&2.20&7.03e-8&2.95\\
		\hline
		1/128&2.68e-3&1.06&4.13e-6&1.98&6.46e-6&2.19&8.85e-9&2.99\\
		\hline\hline
	\end{tabular}\vskip-0.3cm
	\caption{DG and $H^1$ norms of $\widetilde{e}_\psi$ and $\widetilde{e}_\eta$ in Example 3 for $m=4$.}
\end{table}

\begin{table}[H]
	\centering
	\renewcommand\arraystretch{1.0}
	\begin{tabular}{l|l|l|l|l|l|l}
		\hline\hline
		$h$&$\|\widetilde{e}_{\psi}\|_{L^2(\Omega)}$&order&$\|\widetilde{e}_{\eta}\|_{L^2(\Omega)}$&order&$\|\widetilde{e}_{u}\|_{L^2(\Gamma)}$&order\\
		\hline
		1/8&1.43e-4&&1.84e-6&&5.54e-3&\\
		\hline
		1/16&3.32e-5&2.11&1.84e-7&3.32&2.54e-3&1.12\\
		\hline
		1/32&4.43e-6&2.91&1.60e-8&3.52&9.77e-4&1.38\\
		\hline
		1/64&5.47e-7&3.02&1.21e-9&3.73&3.43e-4&1.51\\
		\hline
		1/128&6.30e-8&3.12&9.5e-11&3.67&1.21e-4&1.50\\
		\hline\hline
	\end{tabular}\vskip-0.3cm
	\caption{$L^2$ norms of $\widetilde{e}_\psi,\widetilde{e}_\eta$ and $\widetilde{e}_u$ in Example 3 for $m=4$.}\label{tb18}
\end{table}

So far, we have carried out some numerical experiments for problems with exact manufactured solutions. Since these solutions are smooth, we can observe higher order convergence rates. It is also interesting to consider problems with limited regularity as in Theorem  \ref{thm:regularity} and to verify the convergence rate provided in Theorem \ref{thm3mod}. This can be achieved by considering problems posed on polygonal domains with the largest interior angle greater than 126.283 degrees, as done in this numerical example. From Table \ref{tb13}-\ref{tb18} we can see that, for $h$ small enough, 
the convergence order for the control under $L^2$ norm is nearly $1$ for $m=2$ and is nearly ${3/ 2}$ for $m=3$ and $m=4$, which are consistent with the theoretical results in Theorem \ref{thm3mod}. However, the convergence rates for energy norms and other norms of 
$\widetilde{e}_\psi,\widetilde{e}_\eta$ are higher than  the theoretical results, this is common for Dirichlet boundary control problems (cf. \cite{casas2006dirichlet,may2013dirichlet}), and can be improved by duality arguments (cf. \cite{may2013dirichlet}) for problems without control constraints.  


\section{Conclusion}\label{sec5}
This work applies the reconstructed discontinuous 
approximation method to the biharmonic Dirichlet boundary
control problems with the control imposed on the 
normal derivative value. 
We establish the well-posedness of the discrete 
problem and derive a priori error 
estimates in the energy and $L^2$ norms. 
The numerical results substantiate the theoretical estimates and demonstrate 
the practical effectiveness of this method.

\bibliographystyle{amsplain}
\bibliography{ref}

@article{chowdhury2015priori,
 title={A framework for the error analysis of discontinuous finite element methods for elliptic optimal control problems and applications to {$C^0$} {IP} methods},
 author={Chowdhury, Sudipto and Gudi, Thirupathi and Nandakumaran, A. K.},
 journal={Numer. Funct. Anal. Optim.},
 fjournal= {Numerical Functional Analysis and Optimization},
 volume={36},
 year={2015},
 number={11},
 pages={1388--1419},
 url={https://doi.org/10.1080/01630563.2015.1068182}
}

@article{li2024RDA,
  title={Preconditioned nonsymmetric/symmetric discontinuous {G}alerkin method for elliptic problem with reconstructed discontinuous approximation},
  author={Li, Ruo and Liu, Qicheng and Yang, Fanyi},
  journal={J. Sci. Comput.},
  fjournal={Journal of Scientific Computing},
  volume={100},
  number={88},
  year={2024},
  publisher={Springer},
  url={https://doi.org/10.1007/s10915-024-02639-6}
}

@article{li2024biharm,
  title={An optimal preconditioned discontinuous {G}alerkin method for biharmonic equation with {$C^0$}-reconstructed approximation},
  author={Li, Ruo and Liu, Qicheng and Yang, Fanyi},
  journal={Adv. Comput. Math.},
  fjournal={Advances in Computational Mathematics},
  volume={52},
  number={4},
  year={2026},
  publisher={Springer},
  url={https://doi.org/10.1007/s10444-025-10277-1}
}

@article{li2007posteriori,
  title={A posteriori error estimates of recovery type for distributed convex optimal control problems},
  author={Li, Ruo and Liu, Wenbin and Yan, Ningning},
  journal={J. Sci. Comput.},
  fjournal={Journal of Scientific Computing},
  volume={33},
  number={2},
  pages={155--182},
  year={2007},
  publisher={Springer},
  url={https://doi.org/10.1007/s10915-007-9147-7}
}

@book{lions1971ocpde,
  title={Optimal {C}ontrol of {S}ystems {G}overned by {P}artial {D}ifferential {E}quations},
  author={Lions, Jacques Louis},
  volume={170},
  year={1971},
  publisher={Springer},
  url={https://link.springer.com/book/9783642650260}
}

@book{manzoni2021ocpde,
  title={Optimal {C}ontrol of {P}artial {D}ifferential {E}quations},
  author={Manzoni, Andrea and Quarteroni, Alfio and Salsa, Sandro},
  year={2021},
  publisher={Springer},
  url={https://doi.org/10.1007/978-3-030-77226-0}
}

@article{falk1973ocp,
  title={Approximation of a class of optimal control problems with order of convergence estimates},
  author={Richard S. Falk},
  journal={J. Math. Anal. Appl.},
  fjournal={Journal of Mathematical Analysis and Applications},
  volume={44},
  number={1},
  pages={28--47},
  year={1973},
  url={https://doi.org/10.1016/0022-247X(73)90022-X}
}

@article{geveci1979ocp,
  title={On the approximation of the solution of an optimal control problem governed by an elliptic equation},
  author={Geveci, Tun{\c{c}}},
  journal={RAIRO. Anal. Num{\'e}r.},
  fjournal={RAIRO. Analyse Num{\'e}rique},
  volume={13},
  number={4},
  pages={313--328},
  year={1979},
  publisher={EDP Sciences},
  url={https://www.numdam.org/item/M2AN_1979__13_4_313_0}
}

@article{brenner2019neumann,
  title={{$C^0$} interior penalty methods for an elliptic state-constrained optimal control problem with {N}eumann boundary condition},
  author={Brenner, Susanne C. and Sung, Li-yeng and Zhang, Yi},
  journal={J. Comput. Appl. Math.},
  fjournal={Journal of Computational and Applied Mathematics},
  volume={350},
  pages={212--232},
  year={2019},
  publisher={Elsevier},
  url={https://doi.org/10.1016/j.cam.2018.10.015}
}

@article{casas2012neumann,
  title={Error estimates for the numerical approximation of {N}eumann control problems governed by a class of quasilinear elliptic equations},
  author={Casas, Eduardo and Dhamo, Vili},
  journal={Comput. Optim. Appl.},
  fjournal={Computational Optimization and Applications},
  volume={52},
  number={3},
  pages={719--756},
  year={2012},
  publisher={Springer},
  url={https://doi.org/10.1007/s10589-011-9440-0}
}

@article{casas2008neumann,
  title={Error estimates for the numerical approximation of {N}eumann control problems},
  author={Casas, Eduardo and Mateos, Mariano},
  journal={Comput. Optim. Appl.},
  fjournal={Computational Optimization and Applications},
  volume={39},
  number={3},
  pages={265--295},
  year={2008},
  publisher={Springer},
  url={ https://doi.org/10.1007/s10589-007-9056-6}
}

@article{casas2006dirichlet,
  title={Error estimates for the numerical approximation of {D}irichlet boundary control for semilinear elliptic equations},
  author={Casas, Eduardo and Raymond, Jean-Pierre},
  journal={SIAM J. Control Optim.},
  fjournal={Journal on Control and Optimization},
  volume={45},
  number={5},
  pages={1586--1611},
  year={2006},
  publisher={SIAM},
  url={https://doi.org/10.1137/050626600}
}

@article{may2013dirichlet,
	title={Error analysis for a finite element approximation of elliptic {D}irichlet boundary control problems},
	author={May, Sandra and Rannacher, Rolf and Vexler, Boris},
	journal={SIAM J. Control Optim.},
	fjournal={Journal on Control and Optimization},
	volume={51},
	number={3},
	pages={2585--2611},
	year={2013},
	publisher={SIAM},
	url={https://doi.org/10.1137/080735734}
}

@article{of2015dirichlet,
	title={An energy space finite element approach for elliptic {D}irichlet boundary control problems},
	author={Of, G{\"u}nther and Phan Xuan, Thanh and Steinbach, Olaf},
	journal={Numer. Math.},
	fjournal={Numerische Mathematik},
	volume={129},
	number={4},
	pages={723--748},
	year={2015},
	publisher={SIAM},
	url={https://doi.org/10.1007/s00211-014-0653-x}
}

@article{chowdhury2017dirichlet,
  title={Error bounds for a {D}irichlet boundary control problem based on energy spaces},
  author={Chowdhury, Sudipto and Gudi, Thirupathi and Nandakumaran, A. K.},
  journal={Math. Comp.},
  fjournal={Mathematics of Computation},
  volume={86},
  number={305},
  pages={1103--1126},
  year={2017},
  publisher={AMS},
  url={http://dx.doi.org/10.1090/mcom/3125}
}

@article{gong2022dirichlet,
  title={Analysis and approximations of {D}irichlet boundary control of {S}tokes flows in the energy space},
  author={Gong, Wei and Mateos, Mariano and Singler, John R. and Zhang, Yangwen},
  journal={SIAM J. Numer. Anal.},
  fjournal={SIAM Journal on Numerical Analysis},
  volume={60},
  number={1},
  pages={450--474},
  year={2022},
  publisher={SIAM},
  url={}
}

@article{chen2023dirichlet,
  title={A new global divergence free and pressure-robust {HDG} method for tangential boundary control of {S}tokes equations},
  author={Chen, Gang and Gong, Wei and Mateos, Mariano and Singler, John R. and Zhang, Yangwen},
  journal={Comput. Methods Appl. Mech. Engrg.},
  fjournal={Comput. Methods Appl. Mech. Engrg.},
  volume={409},
  number={},
  pages={No. 115837},
  year={2023},
  publisher={Elsevier},
  url={}
}

@article{georgoulis2009,
  title={Discontinuous {G}alerkin methods for the biharmonic problem},
  author={Georgoulis, Emmanuil H. and Houston, Paul},
  journal={IMA J.
Numer. Anal.},
  fjournal={IMA Journal of
Numerical Analysis},
  volume={29},
  number={3},
  pages={573--594},
  year={2009},
  publisher={Oxford},
  url={}
}

@article{gong2019adaptive,
  title={A convergent adaptive finite element method for elliptic {D}irichlet boundary control problems},
  author={Gong, Wei and Liu, Wenbin and Tan, Zhiyu and Yan, Ningning},
  journal={IMA J. Numer. Anal.},
  fjournal={IMA Journal of Numerical Analysis},
  volume={39},
  number={4},
  pages={1985--2015},
  year={2019},
  publisher={Oxford University Press},
  url={https://doi.org/10.1093/imanum/dry051}
}

@article{hinze2005varydiscre,
  title={A variational discretization concept in control constrained optimization: the linear-quadratic case},
  author={Hinze, Michael},
  journal={Comput. Optim. Appl.},
  fjournal={Computational Optimization and Applications},
  volume={30},
  number={1},
  pages={45--61},
  year={2005},
  publisher={Springer},
  url={https://doi.org/10.1007/s10589-005-4559-5}
}

@book{hinze2008ocpde,
  title={Optimization with {PDE} {C}onstraints},
  author={Hinze, Michael and Pinnau, Ren{\'e} and Ulbrich, Michael and Ulbrich, Stefan},
  year={2008},
  publisher={Springer},
  url={https://doi.org/10.1007/978-1-4020-8839-1}
}

@book{brenner2008book,
  title={The {M}athematical {T}heory of {F}inite {E}lement {M}ethods},
  author={Brenner, Susanne C. and Scott, L. Ridgway},
  year={2008},
  publisher={Springer},
  url={https://doi.org/10.1007/978-0-387-75934-0}
}

@article{chowdhury2017bipriori,
  title={A {$C^0$} interior penalty method for the {D}irichlet control problem governed by biharmonic operator},
  author = {Chowdhury, Sudipto and Gudi, Thirupathi},
  journal={J. Comput. Appl. Math.},
  fjournal={Journal of Computational and Applied Mathematics},
  volume={317},
  pages={290--306},
  year={2017},
  doi={https://doi.org/10.1016/j.cam.2016.12.005}
}

@article{chowdhury2024biprioripos,
  title={Modified {$C^0$} interior penalty analysis for fourth order {D}irichlet boundary control problem and a posteriori error estimates},
  author = {Chowdhury, Sudipto and Garg, Divay and Shokeen, Ravina},
  journal={Math. Comput. Simul.},
  fjournal={Mathematics and Computers in Simulation},
  volume={219},
  pages={185--211},
  year = {2024},
  doi={https://doi.org/10.1016/j.matcom.2023.12.025}
}

@article{garg2025biprioripos,
  title={Discontinuous {G}alerkin method for {D}irichlet boundary control problem governed by biharmonic operator},
  author={Garg, Divay and Porwal, Kamana},
  journal={J. Sci. Comput},
  fjournal={Journal of Scientific Computing},
  volume={103},
  number={23},
  year={2025},
  publisher={Springer},
  doi={https://doi.org/10.1007/s10915-025-02841-0}
}

@article{chowdhury2014bidistributed,
  title={An interior penalty method for distributed optimal control problems governed by the biharmonic operator},
  author={Gudi, Thirupathi and Nataraj, Neela and Porwal, Kamana},
  journal={Comput. Math. Appl.},
  fjournal={Computers \& Mathematics with Applications},
  volume={68},
  number={12, Part B},
  pages={2205--2221},
  year={2014},
  publisher={Elsevier},
  doi={https://doi.org/10.1007/s10092-012-0063-3},
}

@article{frei2013bidistributed,
  title={A priori error estimates for the finite element discretization of optimal distributed control problems governed by the biharmonic operator},
  author={Frei, Stefan and Rannacher, Rolf and Wollner, Winnifried},
  journal={Calcolo},
  fjournal={Calcolo},
  volume={50},
  pages={165--193},
  year={2013},
  publisher={Springer},
  doi={https://doi.org/10.1016/j.camwa.2014.08.012},
}

@book{grisvard2011regupara,
  author={Grisvard, Pierre},
  title={Elliptic {P}roblems in {N}onsmooth {D}omains},
  year={2011},
  publisher={SIAM},
  url={https://epubs.siam.org/doi/abs/10.1137/1.9781611972030}
}

@article{caoyangbiharmonic,
	title={Ciarlet-Raviart mixed finite element approximation for an optimal control problem governed by the first bi-harmonic equation},
	author={Cao, Weidong and Yang, Danping},
	journal={J. Comput. Appl. Math.},
	fjournal={J. Comput. Appl. Math.},
	volume={233},
	pages={372--388},
	year={2009},
	publisher={Elsevier},
	doi={https://doi.org/10.1016/j.cam.2009.07.039},
}

@article{aradacasastroltzsch,
	title={Error estimates for the numerical approximation of a semilinear elliptic control problem},
	author={ Arada, Nadir and Casas, Eduardo and Tr\"oltzsch, Fredi},
	journal={Comput. Optim. Appl.},
	fjournal={Comput. Optim. Appl.},
	volume={23},
	pages={201--229},
	year={2002},
	publisher={Springer},
	doi={},
   }

@article{blumrannacher,
	title={On the boundary value problem of the biharmonic operator on domains with angular corners},
	author={Blum, Heribert and Rannacher, Rolf and Leis, Rolf},
	journal={Math. Methods Appl. Sci.},
	fjournal={Math. Methods Appl. Sci.},
	volume={2},
	pages={556--581},
	year={1980},
	publisher={Teubner},
	doi={https://doi.org/10.1002/mma.1670020416},
}

@article{gudi2010,
  title={A new error analysis for discontinuous finite element methods for linear elliptic problems},
  author={Gudi, Thirupathi},
  journal={Math. Comput.},
  fjournal={Mathematics of Computations},
  volume={79},
  number={272},
  pages={2169--2189},
  year={2010},
  publisher={Americal Mathematical Society},
  url={}
}

@article{hu2020dirichlet,
  title={Analysis of a hybridizable discontinuous {G}alerkin scheme for the tangential control of the {S}tokes system},
  author={Hu, Weiwei and Gong, Wei and Mateos, Mariano and Singler, John R. and Zhang, Yangwen},
  journal={ESAIM Math. Model. Numer. Anal.},
  fjournal={ESAIM Math. Model. Numer. Anal.},
  volume={54},
  number={6},
  pages={2229--2264},
  year={2020},
  publisher={ESAIM},
  url={}
}

\end{document}